\documentclass{amsart}

\usepackage[all]{xy}

\newcommand{\BBox}{\unitlength1.5ex\begin{picture}(1,1)
\put(0,0){\line(1,0){1}}\put(0,0){\line(0,1){1}}\put(1,1){\line(0,-1){1}}
\put(1,1){\line(-1,0){1}}
\end{picture}}

\theoremstyle{plain}
\newtheorem{thm}{Theorem}[section]

\newtheorem{lemma}[thm]{Lemma}

\newtheorem{cor}[thm]{Corollary}
\newtheorem{question}[thm]{Question}

\newtheorem{subthm}{Theorem}[subsection]

\newtheorem{sublem}[subthm]{Lemma}
\newtheorem{subcor}[subthm]{Corollary}

\theoremstyle{definition}
\newtheorem{dfn}[thm]{Definition}
\newtheorem{subdfn}[subthm]{Definition}

\newtheorem{subhypo}[subthm]{Hypothesis}
\newtheorem{subnot}[subthm]{Notation}

\theoremstyle{remark}
\newtheorem{rem}[thm]{Remark}

\newtheorem{subrem}[subthm]{Remark}

\newcommand{\HH}{\mathrm{H}}
\newcommand{\PP}{\mathbb{P}}
\newcommand{\Cl}{\mathrm{Cl}}

\newcommand{\CT}{\mathrm{CT}}

\def\arrow#1#2{\smash{\mathop{\longrightarrow}\limits^{#1}_{#2}}} 

\begin{document}

\title{Galois structure of homogeneous coordinate rings}

\author{Frauke M. Bleher}
\address{F.B.: Department of Mathematics\\University of Iowa\\
Iowa City, IA 52242-1419}
\email{fbleher@math.uiowa.edu}
\author{Ted Chinburg}
\address{T.C.: Department of Mathematics\\University of
Pennsylvania\\Philadelphia, PA
19104-6395}
\email{ted@math.upenn.edu}
\thanks{The first author was supported in part by  
NSF Grant DMS01-39737.
The second author was supported in part by  NSF Grant DMS00-70433.}

\begin{abstract}
Suppose $G$ is a finite group acting on a projective scheme $X$ over a commutative Noetherian ring $R$.  We study the $RG$-modules
$\HH^0(X,\mathcal{F} \otimes \mathcal{L}^n)$ when $n \ge 0$, and $\mathcal{F}$ and $\mathcal{L}$ are coherent $G$-sheaves on $X$ such that $\mathcal{L}$
is an ample line bundle.  We show that the classes of these
modules in the Grothendieck group $G_0(RG)$ of all finitely generated
$RG$-modules lie in a finitely generated subgroup.  Under various
hypotheses, we show that there is a finite set of indecomposable $RG$-modules
such that each $\HH^0(X,\mathcal{F} \otimes \mathcal{L}^n)$ is a direct
sum of these indecomposables, with multiplicites given by
generalized Hilbert polynomials for $n >> 0$.  \end{abstract}

\maketitle


\section{Introduction}
\label{s:intro}
\setcounter{equation}{0}
Let $G$ be a finite group, and suppose that $X$ is a projective scheme over
a commutative Noetherian ring $R$ with an action of $G$ over $R$.  The primary objective
of coherent equivariant Riemann-Roch Theorems is to determine equivariant Euler characteristics associated
to coherent $G$-sheaves $\mathcal{F}$ on $X$.  This problem has motivated
much of the development of 
equivariant $K$-theory and intersection theory, and has a substantial literature
(c.f. \cite{Kock}, \cite{VV}, \cite{VV2}, \cite{Ed}, \cite{Ed2},  
\cite{pappas2}, 
and their references).   One application of
results on Euler characteristics is to study the $RG$-modules  $\HH^0(X,\mathcal{F} \otimes \mathcal{L}^n)$ when $\mathcal{L}$ is an ample $G$-equivariant line bundle on
$X$ and $n >> 0$. The most precise information one would like concerning
$\HH^0(X,\mathcal{F} \otimes \mathcal{L}^n)$ would be an explicit description of this module as a direct sum of 
indecomposable $RG$-modules.   If $RG$ is not semi-simple, this is a more subtle
problem than determining Euler characteristics.  In this paper we will produce such descriptions under various geometric hypotheses on $X$, and we prove some general finiteness results concerning the classes in $G_0(RG)$ defined by
the $\HH^0(X,\mathcal{F} \otimes \mathcal{L}^n)$ as $n$ varies.  The relation of our results to some particular results in the literature will be discussed at the end of this introduction.   

To state one of our main finiteness results, let   $S(X,\mathcal{F},\mathcal{L})$ be the graded $RG$-module \linebreak[4]$\bigoplus_{n\ge 0}
\HH^0(X,\mathcal{F} \otimes \mathcal{L}^n)$.
A graded $RG$-module $M = \bigoplus_{n \in \mathbb{Z}} M_n$ will be said
to be indecomposably finite if there is a finite set $U$ of finitely generated indecomposable
$RG$-modules such that each $M_n$ is isomorphic to the direct sum
of elements of $U$.  A weaker condition is that there is a finitely generated subgroup of the Grothendieck group
$G_0(RG)$ of all finitely generated $RG$-modules which contains the class of $M_n$ for all $n$.
If this is true,  
we will say $M$ is finitely generated in $G_0(RG)$.   

\begin{thm}
\label{thm:thebigzero}
The graded $RG$-module $S(X,\mathcal{F},\mathcal{L})$ is finitely generated
in $G_0(RG)$. 
\end{thm}

It remains to consider when $S(X,\mathcal{F},\mathcal{L})$ is indecomposably finite, and to determine explicit
decompositions of the graded summands into direct sums of
indecomposable $RG$-modules.  In general,  $S(X,\mathcal{F},\mathcal{L})$ is not indecomposably finite even when $R$ is a Dedekind ring 
(see \S\ref{ss:counterex}). A further difficulty is that 
$RG$ need not have the Krull-Schmidt property, so the multiplicities of  indecomposables in a given $RG$-module are not well-defined.
But one can still consider sufficient conditions on $X$, $\mathcal{F}$ and
$\mathcal{L}$ so that the graded  summands of $S(X,\mathcal{F},\mathcal{L})$ can be constructed in the following way.

\begin{dfn}
\label{def:polydef}
A \emph{polynomial description} $\mathfrak{P}$ of a graded $RG$-module $M = \bigoplus_{n \in \mathbb{Z}} M_n$ consists of the following data:
\begin{enumerate}
\item[i.] An integer $m > 0$;
\item[ii.] a finite set $U$  of finitely generated non-zero $RG$-modules; and
\item[iii.]  a polynomial $P_{a,T}(t)   \in \mathbb{Q}[t]$  for each $T \in U$ 
and each integer $a$ in the range $0 \le a < m$ with
$P_{a,T}(t)\neq 0$ for at least one $a$;
\end{enumerate}
satisfying the following conditions:
\begin{enumerate}
\item[a.] All elements of $U$ are indecomposable;
\item[b.] for all sufficiently large $t \in \mathbb{Z}$, $P_{a,T}(t)$ is an integer; and
\item[c.] for all $t >> 0$  and $T \in U$ one has $P_{a,T}(t) \ge 0$, and there is an isomorphism of $RG$-modules
\begin{equation}
\label{eq:sumcond}
M_{tm + a} \cong \bigoplus_{T \in U} T^{P_{a,T}(t)}.
\end{equation}
\end{enumerate}
The degree $d(\mathfrak{P})$ of
this description is the maximum of the degrees of the $P_{a,T}(t)$ in $t$.
\end{dfn}
Note that if finitely generated $RG$-modules have the Krull-Schmidt
property, then the polynomials $P_{a,T}(t)$ are uniquely determined
if they exist.
 
 It is useful to consider the following variant of this definition.  

\begin{dfn}
\label{def:grothdef}
Let $\mathfrak{P}$ be the data specified in Definition
\ref{def:polydef}(i)-(iii).  
Let  $\mathcal{G}$ be a Grothendieck group of finitely generated $RG$-modules, and  suppose the classes 
$[T]$ and $[M_{tm + a}]$ are well-defined in $\mathcal{G}$ for $T \in U$, for all
$a$ and all sufficiently large $t$.  We say that $\mathfrak{P}$ is a 
\emph{$\mathcal{G}$-polynomial description}
of $M$ if instead of (c) in Definition \ref{def:polydef} we have
\begin{equation}
\label{eq:sumclasscond}
[M_{tm + a}] = \sum_{T \in U} P_{a,T}(t)\cdot [T ]\quad \mbox{in $\mathcal{G}$ for $t>>0$}
\end{equation}
where we no longer require $P_{a,T}(t) \ge 0$.  
\end{dfn}

\begin{question} 
\label{qu:whoatethecat}
When does $S(X,\mathcal{F},\mathcal{L})$
have a polynomial  description $($resp. a  $\mathcal{G}$-polynomial de\-scrip\-tion$)$? 
\end{question}

We obtain the following answer to Question \ref{qu:whoatethecat} concerning $\mathcal{G}=G_0(RG)$. In particular, this implies Theorem \ref{thm:thebigzero}.
\begin{thm}
\label{thm:thebigtwo}
There is a $G_0(RG)$-polynomial description of $S(X,\mathcal{F},\mathcal{L})$ of degree bounded by $\mathrm{dim}(\mathrm{supp}(\mathcal{F}))$.
\end{thm}

Our other results concerning Question   \ref{qu:whoatethecat}
have to do with the following two cases: (i) the action of 
$G$ on $X$ is tame, and (ii) $R$ is a field.

We first assume the action of $G$
on $X$ is tame. This means that for each point $x \in X$ the order of 
the inertia group $I_x$ of $x$ is relatively prime to the residue characteristic of $x$.
Let $\CT(RG)$
be the Grothendieck group of all finitely generated $RG$-modules which have
finite projective dimension as $\mathbb{Z}G$-modules.  We allow arbitrary
exact sequences of $RG$-modules in defining the relations in $\CT(RG)$.

\begin{thm}
\label{thm:labelitdummy}  Suppose the action of $G$ on $X$ is tame.
Then $S(X,\mathcal{F},\mathcal{L})$
has a $\CT(RG)$-polynomial description of degree bounded by $\mathrm{dim}(\mathrm{supp}(\mathcal{F}))$.
\end{thm}

\begin{thm}
\label{thm:specialtame}
With the assumptions of Theorem $\ref{thm:labelitdummy}$, 
suppose additionally that $R$ is either a field, or a complete discrete valuation ring, or the ring 
of integers of a number field. Assume 
$\mathcal{F}$ is flat as a sheaf of $R$-modules.  Then 
$S(X,\mathcal{F},\mathcal{L})$ is an indecomposably finite projective $RG$-module having a polynomial description. If $R$ is a field, the degree of this description is $\mathrm{dim}(\mathrm{supp}(\mathcal{F}))$; otherwise the degree  
is $\mathrm{dim}(\mathrm{supp}(\mathcal{F}))-1$.
\end{thm}

We now drop the assumption that the action of $G$ on $R$ is tame, and
assume that $R=k$ is a field of characteristic $p$, with $p=0$ allowed.  
In particular, $X$ is a projective scheme over
$k$ with a $G$-action over $k$, $\mathcal{L}$ is an ample $G$-equivariant
line bundle on $X$, and  $\mathcal{F}$ is a non-zero coherent 
$G$-sheaf on $X$.

\begin{dfn}
\label{def:mostly} 
Let $k$ be a field, and let
$M = \bigoplus_{n \in \mathbb{Z}} M_n$ be a graded $kG$-module.  
Let $P(M_n)$ (resp. $F(M_n)$) be a maximal
projective (resp. free) summand of $M_n$, and define $P(M)=\bigoplus_{n \in \mathbb{Z}} P(M_n)$ (resp. $F(M)=\bigoplus_{n \in \mathbb{Z}} F(M_n)$).
\begin{enumerate}
\item[i.] 
 We will
say $M$ has \emph{polynomial growth} if there is an integer $m > 0$
for which the following is true.  For each  integer $a$ in the range $0 \le a < m$
there is a polynomial $Q_a(t)\in \mathbb{Q}[t]$ such that
$$\mathrm{dim}_k(M_{tm+a}) = Q_a(t)$$
for $t >> 0$. Define $d(M)$ to be the maximum of the degrees of the $Q_a(t)$ in $t$. 
\item[ii.] 
We will say $M$ is \emph{projective
$($resp. free$\,)$ up to terms of degree $c \ge 0$} 
if $\mathrm{dim}_k(M_n/P(M_n))=O(n^c)$
(resp. $\mathrm{dim}_k(M_n/F(M_n))=O(n^c)$).  
\item[iii.] If $P(M_n) = M_n$ (resp. $F(M_n) = M_n$) for $n >> 0$, we will
say $M$ is \emph{eventually projective $($resp. eventually free$\,)$}.
\end{enumerate}
\end{dfn}

Let $G_p = G$ if $p = \mathrm{char}(k)=0$, and otherwise let $G_p$ be a Sylow $p$-subgroup of $G$.
Let $B\subset X$ (resp. $B_p\subset X$) be the ramification locus of the quotient morphism $\pi:X \to Y=X/G$ (resp. of $\pi_p:X \to Y_p=X/G_p$).

\begin{thm}
\label{thm:meow}
Suppose $R=k$ is a field of characteristic $p$, with $p=0$ allowed. 
Define $d=\mathrm{dim}(\mathrm{supp}(\mathcal{F}))$, 
$c = \mathrm{dim}(\mathrm{supp}(\mathcal{F}) \cap B)$, and $c_p = \mathrm{dim}(\mathrm{supp}(\mathcal{F}) \cap B_p)$.
\begin{enumerate}
\item[i.] Then $S(X,\mathcal{F},\mathcal{L}) $ and $P(S(X,\mathcal{F},\mathcal{L}) )$ have polynomial growth, and
$$\mathrm{max}\{c_p,d(P(S(X,\mathcal{F},\mathcal{L}) ))\} = d(S(X,\mathcal{F},\mathcal{L}) ) = d.$$
\item[ii.] If $p=0$ $($resp. $\mathrm{supp}(\mathcal{F}) \cap B_p = \emptyset$$\,)$ 
then $S(X,\mathcal{F},\mathcal{L})$ is projective $($resp. eventually projective$\,)$. 
Otherwise, $S(X,\mathcal{F},\mathcal{L}) $ is projective up to terms of degree $c_p$.
\item[iii.] If $\mathrm{supp}(\mathcal{F}) \cap B = \emptyset$ then 
 $S(X,\mathcal{F},\mathcal{L}) $ is eventually free.  
Otherwise, $S(X,\mathcal{F},\mathcal{L}) $ is free up to terms of degree $c$.
\end{enumerate}
\end{thm}

One rationale for considering how the projectivity or freeness of 
the graded pieces of  $S(X,\mathcal{F},\mathcal{L})$ depends on ramification loci
is that this leads to a proof of the indecomposable finiteness
of $S(X,\mathcal{F},\mathcal{L})$ in some additional cases:

\begin{cor}
\label{cor:ohyeah}
Suppose that the intersection of $\mathrm{supp}(\mathcal{F})$
with the ramification locus $B_p$ is either empty or of dimension $0$.  Then  
$S(X,\mathcal{F},\mathcal{L})$ is indecomposably finite and has a polynomial
description of degree $\mathrm{dim}(\mathrm{supp}(\mathcal{F}))$.
\end{cor}

To get further results about indecomposable finiteness, we need to make some additional assumptions about $X$ and $\mathcal{F}$, respectively.  
We say
$X$ is acyclic if $\HH^i(X,\mathcal{O}_X)=0$ for all $i\ge 1$. For example,
$X=\PP^N_k$ is acyclic for all $N$; also all rational surfaces, Enriques surfaces
and the Godeaux surface are acyclic.

\begin{thm} \label{thm:main}
Suppose $R = k$ is a field of positive characteristic $p$, and $X$
is a smooth projective variety over $k$ having a generically free action of $G$
over $k$. 
Let $r$ be a positive integer, and assume one of the following hypotheses 
is satisfied:
\begin{itemize}
\item[i.] $\mathrm{dim}(X)=1$ and $\mathcal{F}$ is an arbitrary non-zero
coherent $G$-sheaf on $X$, or
\item[ii.] $\mathrm{dim}(X)=2$, $X$ is acyclic, $r$ is sufficiently divisible and
$\mathcal{F}=\mathcal{O}_X$, or
\item[iii.] $X=\mathbb{P}^3_k$, $\mathcal{L}=\mathcal{O}_X(1)$ and
$\mathcal{F}=\mathcal{O}_X$.
\end{itemize}
Then $S(X,\mathcal{F},\mathcal{L}^r)$ is an indecomposably finite $kG$-module which
has a polynomial description of degree 
$\mathrm{dim}(\mathrm{supp}(\mathcal{F}))$. 
\end{thm}

Part (i) of Theorem \ref{thm:main}  is immediate from Corollary \ref{cor:ohyeah}.  
Concerning $r$ in part (ii), we prove that it suffices that $r = m_2^2 (\#G) m_1$
where $m_2$ is the maximal prime to $p$ divisor of $\# G$ and $\mathcal{L}^{m_1}$
is a very ample line bundle on $X$ such that $\HH^i(X,\mathcal{L}^{sm_1}) = 0$
for all integers $i,s > 0$.  In Remark \ref{rem:wedontknow}, we discuss
the problem of improving this value of $r$. We assume $\mathcal{F}=\mathcal{O}_X$
in parts (ii) and (iii) so as to be able to identify global sections with elements in
the function field of $X$.

We now describe the organization of the paper together with the
main ingredients involved in the proofs of the above results.

In \S \ref{s:hom} we  
discuss some preliminaries.
In \S \ref{s:bigone}  we prove
Theorems \ref{thm:thebigtwo} and \ref{thm:labelitdummy} using $\gamma$-filtrations. 
We also prove Theorem \ref{thm:specialtame} using special properties
of complete discrete valuation rings and rings of integers of number fields.
We thank G. Pappas for suggesting the use of $\gamma$-filtrations
in proving Theorems \ref{thm:thebigtwo} and \ref{thm:labelitdummy}.  
This technique for studying equivariant Euler characteristics was introduced in his paper \cite{pappas}.
In \S \ref{s:algclosed} we show that 
the proofs of Theorem \ref{thm:meow}, Corollary \ref{cor:ohyeah} and 
Theorem \ref{thm:main} can be
reduced to the case when the base field $k$ is algebraically closed.
In \S \ref{s:composition} we prove Theorem \ref{thm:meow} and Corollary \ref{cor:ohyeah} using Riemann-Roch results of Fulton \cite{Ful} and of Baum, Fulton and Quart \cite{BFQ}.
In \S \ref{s:projvarfields} we prove parts (ii) and (iii) of Theorem \ref{thm:main} in
the following way. We show
there is a splitting of certain exact cohomology sequences
arising from Koszul resolutions of the structure sheaves of zero-dimensional
subschemes of $X$.  These splittings lead to an inductive expression for
the Krull-Schmidt decomposition of the graded pieces of $S(X,\mathcal{O}_X,\mathcal{L})$.
To construct the splittings, we use an argument going back to Katz and
Grothendieck to show that if $X$ is an acyclic smooth variety over a
field of characteristic $p$, there is a closed point of 
$X$ fixed by a Sylow $p$-subgroup of $G$; 
we thank A. Tamagawa for suggesting the use of this argument.  

We now discuss some particular results in the literature which pertain directly to Question \ref{qu:whoatethecat}.

This paper has been motivated in large part by the following result of Karagueuzian and Symonds.  They show in  \cite{symkar1,symkar2,symkar2.5,sym1,sym2}  that if  
$k$ is algebraic over
$\mathbb{F}_p$ (or $k$ is arbitrary in case $N=2$), then $S(\mathbb{P}^N_k,\mathcal{O}_X,\mathcal{O}_X(1))$ 
has a polynomial description of degree $N$. 
This generalizes the work of many authors \cite{AK,Glover,AKo}. 
The condition that $k$ is algebraic over $\mathbb{F}_p$ is a significant 
restriction, since if $k$ is transcendental over $\mathbb{F}_p$, there may exist 
actions of $G$ which are not defined over any algebraic extension of 
$\mathbb{F}_p$ (see Remark \ref{rem:better}).  Theorem \ref{thm:main}
treats the case of $\mathbb{P}_k^3$ for arbitrary $k$ as well as acyclic
(e.g. rational) surfaces.  

The idea  of considering when $M$ is ``asymptotically mostly
projective" or ``asymptotically mostly free", which is similar to Definition \ref{def:mostly}, 
originates in the work of Symonds in \cite{sym1}.  Theorem \ref{thm:meow} generalizes 
\cite[Thms. 1.1 and 1.2]{sym1} by taking into account $B$ and $B_p$.

 In \cite{sym2}, Symonds  defines
a notion of structure theorem for a graded module which is related to the 
polynomial descriptions appearing in Definition \ref{def:polydef}. 
After we sent him our proof of Theorem \ref{thm:thebigtwo}, he proved in \cite{sym2} a  result of a similar nature using a very different homological method.

We would like to thank C.-L. Chai, G. Pappas, P. Symonds, A. Tamagawa and M. Taylor   for helpful comments.


\section{Preliminaries}
\label{s:hom}
\setcounter{equation}{0}

In this section we will fix some notation and hypotheses,
and we recall some well-known facts about equivariant line bundles,
projective embeddings and homogeneous coordinate rings. 

 Let $R$ be a commutative Noetherian ring,
and let $X$ be a projective scheme over $R$.  We suppose that $G$ is a finite group which acts on $X$ over $R$. We will say the action of $G$ on $X$ is generically free if there is an open dense $G$-stable subset of $X$
on which the action of $G$ is \'etale.  We will say the action of $G$ on $X$ is tame if  for all points 
$x \in X$, the order of the inertia group $I_x \subset G$ of $x$ is relatively prime 
to the residue characteristics of $x$.

Let $\mathcal{L}$
be an ample  line bundle on $X$ having an action of $G$ which is compatible with the action
of $G$ on $\mathcal{O}_X$.  Such $\mathcal{L}$ always exist,
by the following argument. From a projective embedding
$X \to \PP^M_R$ over $R$ one can construct a $G$-equivariant
embedding from $X$ to the product $(\PP^M_R)^{\# G}$ over $R$ of $\# G$ copies
of $\PP^M_R$, with $G$ permuting the factors
of $(\PP^M_R)^{\# G}$.  Following this by a Segre
embedding leads to a projective embedding $\iota:X \to \PP^N_R$ which
is compatible with a linear action of $G$ on $\PP^N_R$.  One can
then take $\mathcal{L} = \iota^* \mathcal{O}_{\PP^N_R}(1)$.

Suppose $\mathcal{F}$ is a non-zero 
coherent $\mathcal{O}_X$-module having a 
compatible action of $G$, and define 
$S(X,\mathcal{F},\mathcal{L}) =  \bigoplus_{n\ge 0}
\HH^0(X,\mathcal{F} \otimes \mathcal{L}^n)$. 

Let $\pi:X\to Y= X/G$ be the quotient morphism.
  By the arguments used to show \cite[Thm. 8.1]{Dol}, it
follows that $Y$ is a quasi-projective scheme over $R$. 
Since $X$ is proper and finite over $Y$, this then implies that $Y$ is 
proper and hence a projective scheme over $R$.

\begin{rem}
\label{rem:faithful}
Without loss of generality, we can assume that the action of $G$ on $X$ is faithful because of the following arguments. Suppose $H\le G$ is the kernel of the action of $G$ on $X$.
Then by inflation from $G/H$ to $G$, any $R(G/H)$-module $M$ is also an $RG$-module.
Moreover, $M$ is simple (resp. indecomposable) as $R(G/H)$-module if and only if it has the same property as $RG$-module.
Inflation from $G/H$ to $G$ also gives a natural homomorphism
$G_0(R(G/H))\to G_0(RG)$.

Suppose now that the action of $G$ on $X$ is tame over $R$.
Then for all $x\in X$, the inertia group $I_x$ contains $H$, and thus $\#H$ is relatively prime
to the characteristic of the residue field $k(x)$. Let $\mathcal{Q}$ be a prime in 
$\mathrm{Spec}(R)$. Since $X$ is proper over $R$, there is a point $x\in X$ that
lies over $\mathcal{Q}$. So $R/\mathcal{Q}$ is a subring of $k(x)$, which implies
that $\#H\not\in \mathcal{Q}$. Hence multiplication by $\#H$ is an isomorphism
of $R_{\mathcal{Q}}$ for all $\mathcal{Q}\in\mathrm{Spec}(R)$, and thus $\#H$ 
is a unit in $R$. 
This means that the inflation from $G/H$ to $G$ of a projective (resp. cohomologically trivial) $R(G/H)$-module
is a projective (resp. cohomologically trivial) $RG$-module.   Thus
we have inflation homomorphisms $K_0(R(G/H)) \to K_0(RG)$ and $\CT(R(G/H)) \to \CT(RG)$.
\end{rem}

\begin{rem}
\label{rem:groth}
Theorem \ref{thm:thebigzero} can also be seen in the following purely algebraic context. 
Suppose $G$ acts $R$-linearly  on the free $R$-module $Rx_0\oplus \cdots \oplus Rx_N$
of rank $N+1$ over $R$.  We can extend this action in a unique way
to an action of $G$ by $R$-algebra automorphisms on the polynomial ring $A = R[x_0,\cdots,x_N]$.  Suppose $I$
is a $G$-stable homogeneous ideal of $A$, 
and let $S = \bigoplus_{n \ge 0} S_n$
be the graded ring $A/I$, where the grading is by degree. Suppose $M$ is a graded $RG$-module 
such that there is an action of $S$ on $M$ which respects the gradings of $S$ and of $M$, and 
$G$ acts on $S$ and $M$ compatibly. Then it follows from Theorem \ref{thm:thebigzero} and  \cite[Thm. 2.3.1]{EGA3} that  $M$ is finitely generated in $G_0(RG)$.
\end{rem}


\section{Projective schemes over arbitrary commutative Noetherian rings}
\label{s:bigone}
\setcounter{equation}{0}

The object of this section is to prove Theorems \ref{thm:thebigtwo}, 
\ref{thm:labelitdummy} and \ref{thm:specialtame}. In \S\ref{ss:counterex}, we give an example in which $R$ is a Dedekind ring and $S(X,\mathcal{F},\mathcal{L})$ is not indecomposably finite.

We assume that the ring $R$ of \S \ref{s:hom} is an arbitrary commutative Noetherian ring and that $G$ acts faithfully on $X$ over $R$.


\subsection{Euler characteristics}
\label{ss:euler}

As usual, $G_0(RG)$ denotes the Grothendieck
group of all finitely generated $RG$-modules.
Define $\mathcal{C}$ to be the category of all finitely generated $RG$-modules
which have finite projective dimension as $\mathbb{Z}G$-modules.  We let
$\CT(RG)$ be the Grothendieck group of $\mathcal{C}$ with
respect to all exact sequences.  

Let $X$ be a projective scheme over $R$ with an action of $G$ over $R$.
Let $G_0(G,X)$ (resp. $K_0(G,X)$) be the Grothendieck group of 
all coherent (resp. coherent locally free) $\mathcal{O}_X$-modules 
$\mathcal{T}$ having a compatible
action of $G$.  We then have a naive Euler characteristic homomorphism
\begin{equation}
\label{eq:naive euler}
\chi^{\mathrm{naive}}:G_0(G,X) \to G_0(RG)
\end{equation}
characterized by the fact that 
$$\chi^{\mathrm{naive}}([\mathcal{T}]) = \sum_i (-1)^i [\HH^i(X,\mathcal{T})].$$
In \cite{CE} and \cite[\S 8]{Duke} it is shown if 
the action of $G$ on $X$ is tame 
and $\mathcal{T}$ is a coherent $G$-sheaf
on $X$, then the hypercohomology complex $\HH^\bullet(X,\mathcal{T})$ is isomorphic in
the derived category of the homotopy category of $RG$-modules to a
bounded complex $P^\bullet = \{P^i\}_i$ of objects of $\mathcal{C}$.  
This leads to the existence of a refined
Euler characteristic homomorphism
\begin{equation}
\label{eq:euler}
\chi:G_0(G,X) \to \CT(RG)
\end{equation}
characterized by the fact that
$$\chi([\mathcal{T}]) = \sum_i (-1)^i [P^i].$$
In particular, the image of $\chi([\mathcal{T}])$ in $G_0(RG)$
is $\chi^{\mathrm{naive}}([\mathcal{T}])$.

\begin{sublem}
\label{lem:stupid}  If $\HH^i(X,\mathcal{T}) = 0$ for $i > 0$ then 
$\HH^0(X,\mathcal{T})$ 
is a finitely generated $RG$-module with class 
$[\HH^0(X,\mathcal{T})] = \chi^{\mathrm{naive}}([\mathcal{T}])$ in $G_0(RG)$.  If in addition, the action of $G$ on $X$ is tame, 
then $\HH^0(X,\mathcal{T})$ is an object of $\mathcal{C}$,
and $[\HH^0(X,\mathcal{T})] = \chi([\mathcal{T}])$ in $\CT(RG)$.
\end{sublem}

As in the statements of Theorems \ref{thm:thebigtwo} and \ref{thm:labelitdummy},
let $\mathcal{L}$ be an ample $G$-equivariant line bundle on $X$, and
let $\mathcal{F}$ be a non-zero coherent $\mathcal{O}_X$-module with a
compatible $G$-action. Let $\pi:X \to X/G = Y$ be the quotient morphism.

\begin{sublem}
\label{lem:pullback}  There is an ample line bundle $\mathcal{L}_Y$
on $Y$ and a $G$-equivariant $\mathcal{O}_X$-module isomorphism
between $\pi^*\mathcal{L}_Y$ and $\mathcal{L}^{(\# G)^2}$.
\end{sublem}

\begin{proof} 
 Because $\pi:X \to Y$ is finite, we can find a finite open affine cover $\{V_i\}_i$ of $Y$
 such that the restriction $\mathcal{L}|_{U_i}$ of $\mathcal{L}$ to  $U_i = \pi^{-1}(V_i)$
 is free of rank one over $\mathcal{O}_{U_i}$ on some generator $f_i$.  For each $g \in G$,
 let $g^*\mathcal{L} = \mathcal{O}_X \otimes_{\mathcal{O}_X,g} \mathcal{L}$ be the pullpack of $\mathcal{L}$ via the automorphism $g:X \to X$.
 Then $g^*\mathcal{L}$ is isomorphic to $\mathcal{L}$ as an $\mathcal{O}_X$-module
 via the action of $g$;  this isomorphism sends  $1 \otimes_g f_i$ to $g(f_i)$.  
 This gives an isomorphism 
 $$\mathcal{L}^{\#G} \cong \prod_{g \in G} g^* \mathcal{L}$$
 sending the $G$-invariant element $m_i = \otimes_{g \in G} g(f_i)$ of $\mathcal{L}^{\# G}$
 to $\otimes_{g \in G} (1 \otimes_g f_i)$.  The $m_i$ produce local generators for a line bundle $\mathcal{H}$ on $Y$ such that $\mathcal{L}^{\# G}$ and 
$\pi^*\mathcal{H}$ are isomorphic as $\mathcal{O}_X$-modules. From this and the action of $G$ we obtain
a cohomology class in $\HH^1(G,\mathrm{Aut}_{\mathcal{O}_X}(\mathcal{L}^{\# G}))$.  Since this group is 
annihilated by $\# G$, we can raise $\mathcal{L}^{\# G}$ to the $\# G^{\mathrm{th}}$ power
to make this cohomology class trivial, leading to a $G$-equivariant isomorphism 
between $\mathcal{L}^{(\# G)^2}$ and $\pi^* \mathcal{L}_Y$ when $\mathcal{L}_Y = \mathcal{H}^{\# G}$. 
To show that $\mathcal{L}_Y$ is ample, it suffices to show that some power of $\mathcal{L}_Y$
is ample.  By replacing $\mathcal{L}$ by a sufficiently high power, we can assume that all
of the $f_i$ above are global sections of $\mathcal{L}$, and that the complement $X_{f_i}$ of 
the zero locus of $f_i$ in $X$ is affine.  Then the $m_i$ are global sections of $\mathcal{L}_Y$.
Since $\pi^{-1}(Y_{m_i}) = \bigcap_{g \in G} X_{g(f_i)}$ is affine and $Y_{m_i} = \pi^{-1}(Y_{m_i})/G$, we conclude each $Y_{m_i}$ is affine, so $\mathcal{L}_Y$ is ample on $Y$
by \cite[p. 155]{hart}.
\end{proof}

 \begin{subrem}
\label{rem:whatever}
We need the following combinatorial result which is a slight variant of
\cite[Lemma 4.8]{pappas}.
Suppose $f:\mathbb{Z}\to A$ is a function from $\mathbb{Z}$ to an abelian group $A$. 
Define the difference operator $\Delta$ 
on functions by $\Delta(f)(n) = f(n+1) - f(n)$.
Then $\Delta^{k} (f) (n) = \sum_{\ell=0}^{k} {k\choose\ell} 
(-1)^\ell f(n+k-\ell)$ for all $k\ge 0$.
If $n_0$ is a fixed integer, then
$(1+\Delta)^{n-n_0}(f)(n_0) = f(n)$ for all $n\ge n_0$. 
Thus we obtain the following formula
for $n \ge n_0$:
\begin{equation}
\label{eq:sonice}
f(n) = \sum_{\ell=0}^\infty {n-n_0 \choose \ell} \cdot \Delta^\ell(f)(n_0).
\end{equation}
\end{subrem}

\medskip

\noindent
\textit{Proof of Theorem $\ref{thm:thebigtwo}$ $($resp. Theorem
$\ref{thm:labelitdummy})$.}
By replacing $X$ by the support of $\mathcal{F}$, we can assume 
$d=\mathrm{dim}(X)=\mathrm{dim}(\mathrm{supp}(\mathcal{F}))$.
Let $A=G_0(RG)$ (resp. $A=\CT(RG)$), and let $\Delta$ be the difference operator 
on functions $f:\mathbb{Z} \to A$ defined by 
$\Delta(f)(n) = f(n+1) - f(n)$.  By Lemma \ref{lem:stupid} and Remark 
\ref{rem:whatever},
it will suffice to show that 
if $0 \le j < (\# G)^2$ and $f_j(n) = \chi^{\mathrm{naive}}([\mathcal{F} \otimes \mathcal{L}^{j + (\# G)^2 n}])$
(resp. $f_j(n) = \chi([\mathcal{F} \otimes \mathcal{L}^{j + (\# G)^2 n}])$)
for all $n \ge 0$, then 
\begin{equation}
\label{eq:deldumb}
\Delta^{d+1} (f_j) (n) = \sum_{\ell = 0}^{d+1}  {d+1\choose\ell} 
(-1)^\ell f_j(n+d+1-\ell)=0 \quad {\rm  for }\quad n >> 0
\end{equation}
when
$d = \mathrm{dim}(Y) = \mathrm{dim}(X)$.  
For $\mathcal{L}_Y$ as in Lemma   \ref{lem:pullback},
the class $([\mathcal{L}_Y] - 1)^{d+1}$ lies in the $(d+1)^{\mathrm{st}}$ term in the $\gamma$-filtration on $K_0(Y) = K_0(\{e\},Y)$,
which is trivial (see \cite[Cor. V.3.10 and Appendix to Chap. V]{fultonlang}).  
Thus in $K_0(G,X)$ we have
\begin{eqnarray*}
0 &=& [\mathcal{F}] \cdot [\mathcal{L}]^{j + (\# G)^2 n}\cdot ([\pi^*\mathcal{L}_Y] - 1)^{d+1}\\
& = &[\mathcal{F}] \cdot [\mathcal{L}]^{j + (\# G)^2 n} \cdot ([\mathcal{L}]^{(\# G)^2} - 1)^{d+1}\\
&=& \sum_{\ell = 0}^{d+1} \left ( \begin{array}{c}d+1\\\ell\end{array}\right)
(-1)^\ell [\mathcal{F}] \cdot [\mathcal{L}]^{j + (\# G)^2(n+d+1 - \ell)}
\end{eqnarray*}
Taking naive Euler characteristics (resp. refined Euler characteristics) 
leads to (\ref{eq:deldumb}) and completes the proof.
\hfill $\BBox$


\subsection{Tame actions}
\label{ss:tame}

In this subsection we prove Theorem \ref{thm:specialtame}. 
We assume the action of $G$ on $X$ over $R$ is faithful and tame. 
We will apply the following well-known result.

\begin{sublem}
\label{lem:cancelresult} Let $\mathcal{B}$ be the set of all finitely
generated projective $RG$-modules, and let $\mathcal{B}_s$ 
be the set of locally free $RG$-modules of constant rank at least $s \ge 0$.
\begin{enumerate}
\item[i.] Suppose $R$ is a field, a discrete valuation ring or the ring
of integers of a number field. Then an object of $\mathcal{C}$ is in
$\mathcal{B}$ if and only if it has no $R$-torsion.
\item[ii.] 
Suppose $R$ is a discrete valuation
ring.   Then two elements of $\mathcal{B}$ are isomorphic if and only
if they become isomorphic after tensoring with the fraction
field of $R$.
\item[iii.]  Suppose $R$ is the ring of integers of a number field.
Then the natural homomorphism $K_0(RG) \to \CT(RG)$ is an isomorphism, and
$\mathcal{B} = \mathcal{B}_0$.  If $P$ and $P'$ are two elements of $\mathcal{B}_2$ with the
same class in $\CT(RG)$ then $P$ and $P'$ are isomorphic.  
The kernel of the rank homomorphism $r:K_0(RG) \to \mathbb{Z}$
is defined to be the class group $\Cl(RG)$ of $RG$;  this is a finite
abelian group 
The homomorphism $\mathbb{Z} \to K_0(RG)$
which sends $1$ to $[RG]$ splits $r$, giving an isomorphism
$$K_0(RG) = \mathbb{Z} \oplus \Cl(RG).$$  
\end{enumerate}
\end{sublem}
\begin{proof}
Part (i) follows using similar arguments as in 
\cite[Proof of Thm. IX.5.7]{corpslocaux}. 
Part (ii) is \cite[Thm. 32.1]{CR}. In part (iii), to prove that the 
natural homomorphism $K_0(RG) \to \CT(RG)$ is an isomorphism, one can use 
\cite[Thm. IX.5.8]{corpslocaux} together with part (i). To show that two elements
of $\mathcal{B}_2$ with the same class in $\CT(RG)$ are isomorphic, one can use 
similar arguments as in \cite[Proof of  Cor. 51.30]{CR}.
The remainder of part (iii) can be found e.g. in \cite[\S 1]{ullom}.
\end{proof}

As in the statement of Theorem \ref{thm:specialtame}, $\mathcal{L}$ is
an ample $G$-equivariant line bundle on $X$, and  $\mathcal{F}$ is a 
non-zero coherent $\mathcal{O}_X$-module which
is flat as a sheaf of $R$-modules and has a compatible action of $G$.
To prove Theorem \ref{thm:specialtame}, we need the following result.

\begin{sublem}
\label{lem:allweneed} 
Let $R$ be a  Dedekind ring, and
suppose that $\mathcal{F}$ is flat as a sheaf of $R$-modules.
\begin{enumerate}
\item[i.] The support $\mathrm{supp}(\mathcal{F})$ of $\mathcal{F}$ is  flat  over $\mathrm{Spec}(R)$.
\item[ii.] Assume $R$ is the ring of integers of a number field.  
The $n^{\mathrm{th}}$ graded piece $\HH^0(X,\mathcal{F} \otimes \mathcal{L}^n)$
of $S(X,\mathcal{F},\mathcal{L})$ is a locally free $RG$-module of rank $z(n)$,
where $z(n) \to \infty $ as $n \to \infty$ unless $\mathrm{supp}(\mathcal{F})$
has dimension $1$. In the latter case, $\mathrm{supp}(\mathcal{F})$ is finite over $\mathrm{Spec}(R)$.
\end{enumerate}
\end{sublem}

\begin{proof}  
Since $R$ is Dedekind, to prove
$\mathrm{supp}(\mathcal{F})$ is flat over $R$, it will suffice to show that $\mathrm{supp}(\mathcal{F})$ is the Zariski closure of
$F \otimes \mathrm{supp}(\mathcal{F})$ in $X$, where $F$ is the fraction field of $R$.  Suppose $x$ is a point in $\mathrm{supp}(\mathcal{F})$ not in the Zariski closure
of $F \otimes \mathrm{supp}(\mathcal{F})$.  Then for each point $x'$ of the general fiber of $X$ which has $x$ in its Zariski closure, the stalk $\mathcal{F}_{x'}$ must be trivial.  It follows that $F \otimes_R \mathcal{F}_x = \{0\}$, so  $\mathcal{F}_x$
has non-trivial $R$-torsion, contradicting the assumption that $\mathcal{F}$ is a
sheaf of flat $R$-modules.  This proves part (i). 

Now suppose that $R$ is the ring of integers of a number field.  Since the stalks of 
$\mathcal{F}$ and $\mathcal{F} \otimes \mathcal{L}^n$ are isomorphic at each point of $X$,
$\mathcal{F} \otimes \mathcal{L}^n$ is a flat sheaf of $R$-modules for all $n$.  Hence $\HH^0(X,\mathcal{F} \otimes \mathcal{L}^n)$
is flat over $R$.  For large $n$, this module is cohomologically trivial for $G$
by Lemma \ref{lem:stupid}.  Hence by Lemma \ref{lem:cancelresult},
$\HH^0(X,\mathcal{F} \otimes \mathcal{L}^n)$ is a locally free $RG$-module for large $n$
of some rank $z(n) \ge 0$.  On replacing $X$ by  $\mathrm{supp}(\mathcal{F})$, we can
suppose $X$ is flat over $R$.  Suppose 
$\mathcal{Q} \in \mathrm{Spec}(R)$, and let $k(\mathcal{Q})$ be the residue field of $\mathcal{Q}$.
Since $\mathcal{F} \otimes \mathcal{L}^n$ is a sheaf of   flat $R$-modules, we have for $n >> 0$ that
$$z(n) =  \frac{1}{\#G} \mathrm{dim}_{k(\mathcal{Q})} \HH^0(k(\mathcal{Q}) \otimes_R X,
k(\mathcal{Q}) \otimes_R \mathcal{F} \otimes \mathcal{L}^n)$$
by considering Hilbert polynomials.  Thus for large $n$, $z(n)$ grows as a polynomial
in $n$ of degree equal to the dimension $\mathrm{dim}(k(\mathcal{Q}) \otimes \mathrm{supp}(\mathcal{F}))$
of the intersection of $\mathrm{supp}(\mathcal{F})$ with the fiber of $X$ over $\mathcal{Q}$.  Thus
the Lemma is proved unless all the fibers of $\mathrm{supp}(\mathcal{F})$ over $R$ have dimension $0$.
Suppose now that the latter condition holds.  Then $\mathrm{supp}(\mathcal{F})$ is quasi-finite and
projective over  $R$ so it is finite over $R$ (c.f. \cite[Ex. III.11.2]{hart}).  
 \end{proof}

\medskip

\noindent
\textit{Proof of Theorem $\ref{thm:specialtame}$.}
By replacing $X$ by the support of $\mathcal{F}$, we can assume, by Lemma
\ref{lem:allweneed}, that $X$ is flat over $R$ and 
$d=\mathrm{dim}(X)=\mathrm{dim}(\mathrm{supp}(\mathcal{F}))$. 
By Lemma \ref{lem:stupid}, $\HH^0(X,\mathcal{F}\otimes\mathcal{L}^n)$
is an object of $\mathcal{C}$ for large $n$; moreover, $\HH^0(X,\mathcal{F}\otimes\mathcal{L}^n)$
has no $R$-torsion since $\mathcal{F}\otimes\mathcal{L}^n$ is a sheaf of
flat $R$-modules for all $n$.  Thus 
by Lemma \ref{lem:cancelresult}, it follows that
$\HH^0(X,\mathcal{F}\otimes\mathcal{L}^n)$
is a projective $RG$-module for $n >> 0$. 

Suppose first that $R=F$ is a field. 
Since the natural map
$\iota: K_0(FG)\to G_0(FG)$
is an injection, there is a $\mathbb{Q}$-linear homomorphism 
$$j:\mathbb{Q}\otimes_{
\mathbb{Z}}G_0(FG)\to \mathbb{Q}\otimes_{\mathbb{Z}}K_0(FG)$$
such that $j\circ\iota$ is the identity map when we view $K_0(FG)$ as a 
subgroup of $\mathbb{Q}\otimes_{\mathbb{Z}}K_0(FG)$. Let $\mathrm{Irr}(
FG)$ be a complete set of representatives for the isomorphism classes of 
simple $FG$-modules. Since $G_0(FG)$ is a free abelian group with
basis $\mathrm{Irr}(FG)$,
Theorem \ref{thm:thebigtwo} implies the following. There exists an integer $m>0$
such that for each $0\le a<m$ and for each
$V\in\mathrm{Irr}(FG)$, there exists a polynomial $Q_{a,\mathcal{F},V}(t)\in 
\mathbb{Q}[t]$ of degree less than or equal to $d=\mathrm{dim}(X)$
such that for large $t\in\mathbb{Z}$
$$[\HH^0(X,\mathcal{F}\otimes\mathcal{L}^{tm+a})] = \sum_{V\in\mathrm{Irr}(FG)} 
Q_{a,\mathcal{F},V}(t)\cdot [V]$$
in $G_0(FG)$. Since the finitely generated $FG$-modules have the
Krull-Schmidt property,
applying the homomorphism $j$ to this equality  then shows
that as  $FG$-modules
$$\HH^0(X,\mathcal{F}\otimes\mathcal{L}^{tm+a}) \cong \bigoplus_{V\in\mathrm{Irr}(FG)} T_V^{
P_{a,\mathcal{F},T_V}(t)}$$
where $T_V$ is the projective envelope of $V$ and $P_{a,\mathcal{F},T_V}(t)$
is a $\mathbb{Q}$-linear combination of the $Q_{a,\mathcal{F},V}(t)$
associated to $V\in\mathrm{Irr}(FG)$. Hence $S(X,\mathcal{F},\mathcal{L})$ is 
an indecomposably finite $FG$-module having a  polynomial description of degree
bounded by $d$. Ignoring the $G$-action, we see from the Hilbert polynomial of
$S(X,\mathcal{F},\mathcal{L})$ that the degree of the description is exactly $d$.

Let now $R$ be a complete discrete valuation ring, and let $F$ be its
fraction field.
Then the finitely generated $RG$-modules have the Krull-Schmidt 
property. By Lemma \ref{lem:cancelresult}, the map 
$\tau: K_0(RG)\to K_0(FG)$ induced by tensoring with $F$ over $R$ is
an injection. Thus $1\otimes \tau:\mathbb{Q}\otimes_{\mathbb{Z}}
K_0(RG)\to \mathbb{Q}\otimes_{\mathbb{Z}} K_0(FG)$ has a $\mathbb{Q}$-linear section.
Hence by the first paragraph of the proof,
$S(X,\mathcal{F},\mathcal{L})$ is an indecomposably finite
$RG$-module having a  polynomial description 
if and only if the same is true for
$F \otimes_R S(X,\mathcal{F},\mathcal{L}) = S(F\otimes_R X, F\otimes_R\mathcal{F},F\otimes_R \mathcal{L})$
as an $FG$-module.
The latter fact follows from the second paragraph of the proof.
Since the degree of the polynomial description of $F \otimes_R 
S(X,\mathcal{F},\mathcal{L})$ is $d-1=\mathrm{dim}(F\otimes_R X)$, we
see that this is also the degree of the polynomial description of 
$S(X,\mathcal{F},\mathcal{L})$.

Finally, let $R$ be the ring of integers of a number field with fraction field $F$.  
By Theorem \ref{thm:labelitdummy}, $S(X,\mathcal{F},\mathcal{L})$
has a $\CT(RG)$-polynomial description.
By Lemma \ref{lem:cancelresult}, the quotient $\Cl(RG)$ of
$\CT(RG) = K_0(RG)$ by the subgroup generated by the class of the free module $RG$
is a finite abelian group.  Since polynomials with integer values are $\ell$-adically
continuous functions for all rational primes $\ell$,  the resulting $\Cl(RG)$-polynomial description
of  $S(X,\mathcal{F},\mathcal{L})$ has the following property.  There is an integer $m_0 > 0$
such that the class of the $n^{\mathrm{th}}$ graded piece $\HH^0(X,\mathcal{F} \otimes \mathcal{L}^n)$ of $S(X,\mathcal{F},\mathcal{L})$
in $\Cl(RG)$ depends only on $n \mod m_0$ for $n >> 0$.  Fix an integer $0 \le a < m_0$.
We conclude that there is a projective $RG$-module $Q_a$ such that for $n >> 0$ and $n\equiv a\mod m_0$, the class of $\HH^0(X,\mathcal{F} \otimes \mathcal{L}^{n})$ in 
$\CT(RG) = K_0(RG)$ equals $[Q_a] + (z(n) - r_a) [RG]$ where $z(n)$ (resp. $r_a$)
is the rank of $\HH^0(X,\mathcal{F} \otimes \mathcal{L}^{n})$ (resp. $Q_a$) as a locally free $RG$-module. Here $z(n)=\displaystyle \frac{1}{\#G} \mathrm{dim}_F\,
\HH^0(F\otimes_R X, (F\otimes_R\mathcal{F})\otimes (F\otimes_R\mathcal{L}^n))$.
If $z(n) \to \infty$ as $n \to \infty$, then by Lemma \ref{lem:cancelresult} we conclude that 
$\HH^0(X,\mathcal{F} \otimes \mathcal{L}^{n})$ is isomorphic to $Q_a \oplus  (RG)^{z(n) - r_a}$
as $RG$-modules for $n\equiv a\mod m_0$ and $n >> 0 $. The rank function
$z(n)$ is a polynomial in $n$ of degree equal to $d-1=\mathrm{dim}(
\mathrm{supp}(F\otimes_R\mathcal{F})=\mathrm{dim}(F\otimes_R X)$ 
for $n >> 0$, so we are done in this case.  In the remaining
case,  the support of $\mathcal{F}$ is finite over $R$ by Lemma \ref{lem:allweneed}.  Thus, since we assumed $X=\mathrm{supp}(\mathcal{F})$,
$X$ is flat finite and projective over $R$.
Then  $\mathcal{L}$ defines
an element of the group $\mathrm{Pic}^G(X)$ of all isomorphism classes
of $G$-equivariant line bundles
on $X$, where the group structure of $\mathrm{Pic}^G(X)$ comes from the tensor
product over $\mathcal{O}_X$.  To complete the proof, it will suffice to show $\mathrm{Pic}^G(X)$
is finite, since then the isomorphism class of $\HH^0(X,\mathcal{F}\otimes\mathcal{L}^n)$
as an $RG$-module is periodic as a function of $n$.  
Since $\mathrm{dim}(X) = 1$ and $X$ is flat over $R$, $X = \mathrm{Spec}(A)$ for some $R$-order $A$ in a finite dimensional algebra $N$ over $F$.  Then $\mathrm{Pic}(X) = \mathrm{Pic}(A)$ is finite since $N$ is finite 
dimensional over $\mathbb{Q}$.
There is an injective homomorphism from the kernel of the forgetful map
$\mathrm{Pic}^G(X) \to \mathrm{Pic}(X)$ to $\HH^1(G,\mathrm{Aut}_R(\mathcal{O}_X))$.
Since $\mathrm{Aut}_R(\mathcal{O}_X) = A^*$ is a finitely generated abelian group,
and $G$ is finite, $\HH^1(G,\mathrm{Aut}_R(\mathcal{O}_X))$ must be finite,
so $\mathrm{Pic}^G(X)$ is as well.    
\hfill $\BBox$


\subsection{A counterexample to indecomposable finiteness}
\label{ss:counterex}

Suppose $E$ is an elliptic curve over the complex numbers $\mathbb{C}$.
Let $R$ be the affine ring of the complement in $E$ of the origin.
Then $R$ is a Dedekind ring.   Set $X = \mathrm{Spec}(R)$, and let
$G$ be the trivial group.  Then
$\mathrm{Pic}(X)$ is isomorphic to the group of divisors of degree $0$
on $E$, which is isomorphic to the group of complex points $E(\mathbb{C})$. 
Let $\mathcal{L}$ be a line bundle on $X$ corresponding  to a point
of infinite order of $E(\mathbb{C})$.  Then $I = \HH^0(X,\mathcal{L})$ is a locally
free rank one $R$-module, and $S(X,\mathcal{O}_X,\mathcal{L}) = \bigoplus_{n \ge 0} \HH^0(X,\mathcal{L}^n) = \bigoplus_{n \ge 0} I^n$ is not an indecomposably 
finite $R = RG$-module. Here $I^n$ denotes the $n^{\mathrm{th}}$ power of
the ideal $I$.


\section{Projective schemes over fields: reduction to algebraically closed fields}
\label{s:algclosed}
\setcounter{equation}{0}

Throughout this section we will assume the notations of \S \ref{s:hom}.
Moreover, we will assume that the ring $R$ of \S \ref{s:hom} is a field $k$ and that $G$ acts faithfully on $X$ over $k$.
Let $G_0^{\oplus}(kG)$ be the Grothendieck group  of all finitely
generated $kG$-modules with respect to direct sums.  Since $kG$ has
the Krull-Schmidt property, two finitely generated $kG$-modules are isomorphic if and only
if they have the same class in $G_0^{\oplus}(kG)$.

Define
$P_{k,G}$ (resp. $N_{k,G}$) to be the subgroup of $G_0^{\oplus}(kG)$
generated by projective indecomposable (resp. non-projective
indecomposable) $kG$-modules.  Each of $G_0^{\oplus}(kG)$, $P_{k,G}$ and $N_{k,G}$ are
free abelian groups, and $G_0^{\oplus}(kG)= P_{k,G} \oplus N_{k,G}$.   
The following Lemma is proved using the Krull-Schmidt Theorem.

\begin{lemma}
\label{lem:basext}
Let $k'$ be an algebraic extension of $k$.  
\begin{enumerate}
\item[i.]
The homomorphism $\psi:G_0^{\oplus}(kG) \to G_0^{\oplus}(k'G)$
induced by tensoring with $k'$ over $k$ is injective, and sends 
$P_{k,G}$ into $P_{k',G}$ and $N_{k,G}$ into $N_{k',G}$.  
There is a $\mathbb{Q}$-linear
homomorphism $j:\mathbb{Q} \otimes_{\mathbb{Z}} G_0^{\oplus}(k'G) \to \mathbb{Q} \otimes_{\mathbb{Z}} G_0^{\oplus}(kG)$ such that $j\circ \psi$ is the identity map.
\item[ii.] Suppose $D$ is a finitely generated $kG$-module. If 
$c(k,D)$ is the codimension in $D$ of a maximal free $kG$-module
summand, then $c(k',k'\otimes_k D)=c(k,D)$.
\end{enumerate}
\end{lemma}

\begin{cor}
\label{cor:algext}
A graded $kG$-module $M = \bigoplus_{n \ge 0}M_n$ is indecomposably finite 
$($resp. has a polynomial description of degree $d$, 
resp. is projective $($resp. free$\,)$ up to terms of degree $c$, resp. is eventually
projective $($resp. eventually free$\,)${}$\,)$ if 
the $k'G$-module $k'\otimes_k M$ has the corresponding property 
for some algebraic extension $k'$ of $k$.
\end{cor}

\begin{proof} We suppose first that $k' \otimes_k M$ has a polynomial description
of degree $d$,
and will show $M$ has such a description;  indecomposable finiteness can be treated
in a similar way.  There is then an integer $m > 0$ and a finite set $U'$ of finitely 
generated indecomposable $k'G$-modules $T'$ together with polynomials $P_{a,T'}(t) \in \mathbb{Q}[t]$ such that
for large $t \in \mathbb{Z}$, one has $P_{a,T'}(t) \ge 0$ and
$$k' \otimes_k M_{tm+a} \cong \bigoplus_{T' \in U'} T'^{P_{a,T'}(t)}$$
as $k'G$-modules. Furthermore, the 
degrees of the $P_{a,T'}(t)$ in $t$ are less than or equal to $d$ with equality for at least
one $P_{a,T'}(t)$.
Let $U$ be the finite set of finitely generated indecomposable $kG$-modules $T$ which
occur with non-zero coefficient in $j(T') \in \mathbb{Q} \otimes_{\mathbb{Z}} G_0^{\oplus}(kG)$
for some $T' \in U'$,
where $j: \mathbb{Q} \otimes_{\mathbb{Z}} G_0^{\oplus}(k'G) \to \mathbb{Q} \otimes_{\mathbb{Z}} G_0^{\oplus}(kG)$ is the homomorphism in Lemma \ref{lem:basext}(i).  This Lemma
then shows that as a $kG$-module, 
$$M_{tm+a} \cong \bigoplus_{T \in U} T^{P_{a,T}(t)}$$
where $P_{a,T}(t)$ is a finite $\mathbb{Q}$-linear combination of the $P_{a,T'}(t)$
associated to $T' \in U'$.  Thus 
we arrive at a polynomial description of $M$ of degree $d$.

If $k'\otimes_k M$ is projective (resp. free) up to terms of degree $c$
(resp. eventually projective (resp. eventually free)), then it follows from Lemma
\ref{lem:basext} that the same is true for $M$.
\end{proof}

\begin{lemma}
\label{lem:biggerfield}
To prove Theorems $\ref{thm:meow}$, Corollary $\ref{cor:ohyeah}$ and  Theorem $\ref{thm:main}$, it suffices to consider
the case in
which $k$ is algebraically closed.
\end{lemma}

\begin{proof}
Let $k'$ be an algebraic closure of $k$, and let 
$X' = k' \otimes_k X$ and $\mathcal{L}' = k' \otimes_k \mathcal{L}$.

To prove the statement  Theorem \ref{thm:meow} (resp. Corollary \ref{cor:ohyeah}), 
let
$\mathcal{F}' = k' \otimes_k \mathcal{F}$.
By flat base change, we have
$k'\otimes_k \HH^0(X,\mathcal{F}\otimes\mathcal{L}^n)
= \HH^0(X',\mathcal{F}'\otimes{\mathcal{L}'}^n)$
as $k'G$-modules for all $n$. 
Hence by Corollary \ref{cor:algext}, if Theorem \ref{thm:meow} 
(resp. Corollary \ref{cor:ohyeah}) holds for
$X'$, $\mathcal{L}'$ and $\mathcal{F}'$, then it holds for $X$, $\mathcal{L}$
and $\mathcal{F}$. 

Suppose now that $X$ is a smooth projective variety
over $k$ as in Theorem \ref{thm:main}.  By flat base change, we have
${k'}\otimes_k S(X,\mathcal{L}) = S(X',\mathcal{L}')$
as graded $k'G$-modules.   
Let $\{X_i\}$ be a finite set of representatives for the $G$-orbits of 
irreducible  components of $X'$. Since $X'$ is smooth over $k'$, the $X_i$ are smooth projective varieties over $k'$, and
$X'$ is the disjoint union of the $G$-translates of the $X_i$.  
For each $i$, let $H_i$ be the stabilizer of $X_i$ in $G$,
and let $\mathcal{L}_i$ be the restriction of $ \mathcal{L}'$ to 
$X_i$. Then 
$${k'}\otimes_k S(X,\mathcal{L}) = S(
X',\mathcal{L}')\cong\bigoplus_{i}\mathrm{Ind}_{H_i}^G 
S(X_i,\mathcal{L}_i).$$
Thus, by  Corollary \ref{cor:algext}, to prove Theorem \ref{thm:main}, we can reduce to the case in which $k = k'$.
\end{proof}


\section{Projective schemes over fields: results up to high codimension}
\label{s:composition}
\setcounter{equation}{0}

In this section we will prove Theorem \ref{thm:meow} and Corollary \ref{cor:ohyeah}.
Throughout this section we will assume the notations of \S \ref{s:hom}. Moreover,
we will assume that the ring $R$ of \S \ref{s:hom} is a field $k$. By Lemma
\ref{lem:biggerfield}, we can assume $k$ is algebraically closed.

If $T$ is a finitely generated $kG$-module, we denote by $P(T)$ (resp. $F(T)$)
a maximal projective (resp. free) summand of $T$ and by $P'(T)$ (resp. $F'(T)$) a complement
for $P(T)$  (resp. $F(T)$) in $T$.  The isomorphism
classes of these modules are determined by $T$.
If $M = \bigoplus_{n \in \mathbb{Z}} M_n$ is a graded $kG$-module,
we define
$\mathcal{Z}(M) = \bigoplus_{n \in \mathbb{Z}} \mathcal{Z}(M_n)$ for $\mathcal{Z}\in\{P,P',F,F'\}$.

Let $p$ be the characteristic of $k$, with $p = 0$  allowed.
An element $g \in G$ is $p$-regular if either $p = 0$, or $p > 0$ and $g$ 
has order prime to $p$.  If $p = 0$ we let $\Lambda = k$;  otherwise let $\Lambda$ be  
a complete discrete valuation ring of characteristic $0$ with residue
field $k$ and fraction field sufficiently large relative to $G$.
For each finitely generated $kG$-module
$T$, let $\Phi(T)$ be the ordinary character of $T$ if $p = 0$ and the Brauer character 
of $T$ if $p > 0$. Thus $\Phi(T)$ is a function from the
$p$-regular elements of $G$ to $\Lambda$. Let $F$ be the fraction field of $\Lambda$.
Define $G_p = G$ if $p = 0$, and otherwise let $G_p$ be a Sylow $p$-subgroup
of $G$.


\subsection{Module theoretic results}
\label{ss:meowmodule}
 
\begin{sublem}
\label{lem:projtest}
Suppose $\mathrm{char}(k)=p>0$, and
let $T$ is a finitely generated $kG$-module.  A set
$\{z_i\}_i$ of elements of $T$ is a basis for a free $kG_p$-module
summand of $T$ if and only if the elements $\{\mathrm{Tr}_{G_p} (z_i)\}_i$
are linearly independent over $k$.  In particular,
\begin{equation}
\label{eq:trdim}
\mathrm{dim}_k (P(T)) = \# G_p \cdot \mathrm{dim}_k (\mathrm{Tr}_{G_p} (T))
\end{equation}
\end{sublem}

\begin{proof} Suppose there is a relation $\sum_i \beta_i z_i = 0$ in which
the $\beta_i$ are elements of $kG_p$ which are not all $0$.  We can reduce
to the case in which the radical of $kG_p$ annihilates all of the $\beta_i$.
But then the $\beta_i$ are elements of $k \mathrm{Tr}_{G_p} $ since $G_p$
is a $p$-group, contradicting the assumption that the elements $\{\mathrm{Tr}_{G_p} (z_i)\}_i$
are linearly independent over $k$.
\end{proof}

\begin{sublem}
\label{lem:nolabel}
Suppose $M = \oplus_{n \in \mathbb{Z}} M_n$ is a graded  $kG$-module
in which each $M_n$ is finitely generated.    Suppose there are integer constants  $0 \le \ell_p \le \ell$ such that 
\begin{enumerate}
\item[i.] $\mathrm{dim}_k(P'(M_n))   = O(n^{\ell_p})$; and
\item[ii.]  for each $p$-regular element $g \in G$ which is not the identity element,
there is a polynomial $\rho_g(t) \in F[t]$ of degree less than or equal to 
$\ell$ in $t$ such that $\Phi(M_n)(g) = \rho_g(n)$ for all sufficiently
large $n$.
\end{enumerate}
Then $\mathrm{dim}_k(F'(M_n)) = O(n^\ell)$.  
\end{sublem}

\begin{proof}
Let $V$ be a simple $kG$-module.  If $p = 0$, let $P_V=V$ and let 
$\Xi_V$ be the character of $V$.  If $p > 0$, let $X_V$ be a projective
$\Lambda G$-module whose reduction modulo the maximal ideal of $\Lambda$ is 
the projective $kG$-cover $P_V$ of $V$, 
and let $\Xi_V$ be the character of the $FG$-module $F\otimes_{\Lambda} X_V$.  
For $p \ge 0$, define $G_{p'}$ to be the set of all $p$-regular 
elements of $G$.    
By \cite[\S 2.3 and \S 18.3]{serre}, the multiplicity $m_{n}(V)$ of $V$ in a
composition series of $M_n$ is
\begin{equation}
\label{eq:Vmult}
 \frac{1}{\#G}\;
\sum_{g\in G_{p'}} \Xi_V(g^{-1})\,\Phi(M_n)(g) = 
\frac{1}{\#G}  \Xi_V(e)\cdot \mathrm{dim}_k(M_n) + O(n^\ell)
\end{equation}
Here $e$ is the identity element of $G$, and the equality follows from assumption (ii).  Because of assumption (i), we see that this implies the multiplicity $m'_{n}(V)$ of $V$ 
in a composition series for $P(M_n)$ has the form
\begin{equation}
\label{eq:Vmult2}
m'_n(V) = \frac{1}{\#G}  \Xi_V(e)\cdot \mathrm{dim}_k(P(M_n)) +O(n^\ell).
\end{equation}
Let $P(M_n)=\oplus_V P_V^{c_V(n)}$ where $V$ runs over the isomorphism classes of simple $kG$-modules.  Let $r_n$ be the largest integer less than
or equal to $\frac{1}{\#G} \mathrm{dim}_k(P(M_n))$. Then
$(kG)^{r_n} = \oplus_V P_V^{r_n \mathrm{dim}_k(V)}$. 
In $K_0(kG)$ we thus have
$$[P(M_n)] - r_n[kG] = \sum_V (c_V(n) - r_n \mathrm{dim}_k(V)) [P_V].$$
Since $\Xi_V(e)$ is the multiplicity of $V$ in a composition series of $kG$, the image of this class in $G_0(kG)$ is
$$\sum_V \left(m'_n(V) - \Xi_V(e) r_n\right) [V]$$
where $m'_n(V) - \Xi_V(e) r_n = O(n^\ell)$.  Because the homomorphism $K_0(kG) \to G_0(kG)$
is an injection of finitely generated free abelian groups, this implies
$$c_V(n) - r_n \mathrm{dim}_k(V) =O(n^\ell).$$
We can therefore find an integer $r'_n \ge 0$
such that $0 \le c_V(n) - r'_n \mathrm{dim}_k(V) = O(n^\ell)$.
 Let $F_n$ be a free $kG$-module of rank $r'_n$.  It follows that $P(M_n) 
= F_n \oplus F'_n$
with $\mathrm{dim}_k(F'_n) = O(n^\ell)$.  Because of assumption (i) of 
the Lemma, this implies
$\mathrm{dim}_k(F'(M_n)) = O(n^\ell)$.  
\end{proof}


\subsection{Proof of parts (i) and (ii) of Theorem \ref{thm:meow}}
\label{ss:meow1and2}

As in the statement of Theorem \ref{thm:meow}, $X$ is a projective scheme over 
$k$, $\mathcal{L}$ is an ample $G$-equivariant line bundle on $X$, and 
$\mathcal{F}$ is a non-zero coherent $G$-sheaf on $X$ with 
$d=\mathrm{dim}(\mathrm{supp}(\mathcal{F}))$.  Let $\pi:X\to Y=X/G$,
$\pi_p:X\to Y_p=X/G_p$ and $\zeta:Y_p=X/G_p\to Y=X/G$ be the natural
quotient morphisms.
By Lemma \ref{lem:pullback},  there is an ample line bundle $\mathcal{L}_Y$
on $Y = X/G$ and a $G$-equivariant $\mathcal{O}_X$-module isomorphism
between $\pi^*\mathcal{L}_Y$ and $\mathcal{L}^{m}$ where $m = (\# G)^2$.   
By replacing $m$ by an appropriate multiple, we may assume that
$\mathcal{L}_Y$, $\zeta^*\mathcal{L}_Y$ and $\mathcal{L}^m\cong
\pi^*\mathcal{L}_Y$ are very ample line bundles over $k$.

The fact that
$S(X,\mathcal{F},\mathcal{L})$ has polynomial growth with $d(S(X,\mathcal{F},\mathcal{L}))=d$ follows
from the fact that the Hilbert polynomial of $\mathcal{F}\otimes\mathcal{L}^a$ 
has degree $d$ for every integer $a$. If $\mathrm{char}(k)=p=0$ then
$P(S(X,\mathcal{F},\mathcal{L}))=S(X,\mathcal{F},\mathcal{L})$, and we
are done in this case.

We may replace $X$ by the support of $\mathcal{F}$, so as to be able to
assume $X=\mathrm{supp}(\mathcal{F})$ and $d=\mathrm{dim}(\mathrm{supp}(\mathcal{F}))=\mathrm{dim}(X)$.
Define  $B_p$ to be the ramification locus of the cover $\pi_p$, and let
$c_p = \mathrm{dim}(\mathrm{supp}(\mathcal{F}) \cap B_p) = \mathrm{dim}(\mathrm{supp}( B_p))$.  
If $B_p$ is empty, then $\pi:X \to Y = X/G$ is tamely ramified and $S(X,\mathcal{F},\mathcal{L})$
is eventually projective by Lemmas \ref{lem:stupid} and \ref{lem:cancelresult}(i).  
In particular, $P(S(X,\mathcal{F},\mathcal{L}))$ has polynomial growth and
$d(P(S(X,\mathcal{F},\mathcal{L})))=d$.

In the remainder
of this subsection we suppose $\mathrm{char}(k)=p>0$ and $B_p \ne \emptyset$.  
We will show that  
$P(S(X,\mathcal{F},\mathcal{L}))$
has polynomial growth and that 
$S(X,\mathcal{F},\mathcal{L})$  is projective up to terms of degree $c_p$, 
in the sense that
\begin{equation}
\label{eq:nonoproj} 
\mathrm{dim}_k(P'(\HH^0(X,\mathcal{F}_n)))  = O(n^{c_p}).
\end{equation}
when $\mathcal{F}_n = \mathcal{F}  \otimes \mathcal{L}^n$.  In view of Lemma \ref{lem:projtest}, 
(\ref{eq:nonoproj}) is equivalent to 
\begin{equation}
\label{eq:littleproj}
\mathrm{dim}_k (\HH^0(X,\mathcal{F}_n))  =
\# G_p \cdot \mathrm{dim}_k ( \mathrm{Tr}_{G_p} \HH^0(X,\mathcal{F}_n))  +
O(n^{c_p}).
\end{equation}

 The trace operator $\mathrm{Tr}_{G_p}$ gives an exact sequence
\begin{equation}
\label{eq:Kseq}
0 \to \mathcal{K}_n \to (\pi_{p})_* (\mathcal{F}_n)  \to \mathrm{Tr}_{G_p} (\pi_{p})_* (\mathcal{F}_n) \to 0
\end{equation}
of sheaves on $Y_p$, where  $\mathcal{K}_n$ is coherent.  
Write $n = a + rm$ for some $0 \le a < m$ and $r \ge 0$.  Since
$\mathcal{L}^m$ is $G$-isomorphic to $\pi^* \mathcal{L}_Y$  we find from the projection formula that   $\mathcal{K}_{a + rm} = \mathcal{K}_a \otimes_{\mathcal{O}_{Y_p}} (\zeta^*\mathcal{L}_Y)^r$.
Hence $\HH^1(Y_p,\mathcal{K}_{a + rm}) = 0$ if
$r$ is sufficiently large.  It follows from the cohomology of (\ref{eq:Kseq}) that 
\begin{equation}
\label{eq:asymptot}
\mathrm{Tr}_{G_p} (\HH^0(X,\mathcal{F}_n)) = 
\mathrm{Tr}_{G_p} (\HH^0(Y_p,(\pi_{p})_* (\mathcal{F}_n))) =
\HH^0(Y_p,\mathrm{Tr}_{G_p} (\pi_{p})_* (\mathcal{F}_n)) \quad \mathrm{for}\quad n >> 0.
\end{equation}

 The cover
$U = X - B_p \to V = (X - B_p)/G_p$ is an \'etale  $G_p$-cover.  
This
implies the natural morphism $\pi_p^* \mathrm{Tr}_{G_p} (\pi_p)_* \mathcal{F}_a \to \mathcal{F}_a$
of sheaves on $X$ has kernel and cokernel supported on $B_p$.  This morphism is tensored by $\mathcal{L}^{rm}$ when
$a$ is replaced by $n = a + rm$. Since
$\mathrm{dim}(B_p) = c_p$, we see by taking Euler characteristics and Hilbert polynomials
that
\begin{equation}
\label{eq:zapzap2}
f_{1,a}(r) = \mathrm{dim}_k(\HH^0(X,\mathcal{F}_{a + rm}))  -
\mathrm{dim}_k (\HH^0(X,\pi_p^* \mathrm{Tr}_{G_p} (\pi_p)_* \mathcal{F}_{a + rm}))  
\end{equation}
is a polynomial in $r$ of degree $\le c_p$ for $r >> 0$.  

In view of (\ref{eq:asymptot}) and (\ref{eq:zapzap2}), to show (\ref{eq:littleproj}) we
are reduced to showing that 
\begin{equation}
\label{eq:zapzap}
f_{2,a}(r) = \mathrm{dim}_k (\HH^0(X,\pi_p^* \mathcal{G}_{a + rm} ))  -
\# G_p \cdot \mathrm{dim}_k (\HH^0(Y_p,\mathcal{G}_{a + rm}))  
\end{equation}
is a polynomial of degree $\le c_p$ in $r$ for large $r$ 
when  $\mathcal{G}_n = \mathrm{Tr}_{G_p} (\pi_p)_*\mathcal{F}_n$.  
This will then also prove that $P(S(X,\mathcal{F},\mathcal{L}))$ has polynomial
growth.

We can now use the Riemann-Roch result of \cite[Thm. 18.2]{Ful}
for quasi-projective schemes $Z$ over $k$.   Let 
$$\tau_Z: G_0(Z) \to A_* (Z)_{\mathbb{Q}}$$
be the fundamental homomorphism defined in this Theorem from
the Grothendieck group $G_0(Z)$ of all coherent sheaves on $Z$ to
the rational Chow ring $A_* (Z)_{\mathbb{Q}}$ of $Z$. (In \cite{Ful},
$G_0(Z)$ is denoted by $K_{o}(Z)$.)  Define $b_p = B/G_p$, so
that $V = U/G_p = Y_p - b_p$. From \cite[Prop. 1.8]{Ful}, the closed 
immersion $i:b_p \to Y_p$
and open immersion $j:V \to Y_p$ 
give an exact sequence
\begin{equation}
\label{eq:Az}
A_* ({b_p})_{\mathbb{Q}} \arrow{i_*}{} A_* (Y_p)_{\mathbb{Q}} \arrow{j^*}{} A_* (V)_{\mathbb{Q}}\to 0.
\end{equation}
Let $\pi_p$ stand for each of the quotient morphisms $X \to Y_p = X/G_p$,
$U \to U/G_p = V$ and $B_p \to B/G_p = b_p$.   
Since $U \to V$ is
an \'etale $G_p$-cover, it is smooth of relative dimension $0$, and
\cite[Thm. 18.2]{Ful} shows that
$$\tau_V((\pi_p)_*  \pi_p^* j^* [\mathcal{G}_n]) = 
(\pi_p)_*  \pi_p^*  \tau_V( j^* [\mathcal{G}_n]) = \# G_p \cdot
\tau_V(j^* [\mathcal{G}_n]) \quad \mathrm{in}\quad A_* (V)_{\mathbb{Q}}.$$
Hence (\ref{eq:Az}) shows that 
\begin{equation}
\label{eq:pushit}
 \tau_{Y_p}((\pi_p)_* \pi_p^*[ \mathcal{G}_n]) - \# G_p \cdot
\tau_{Y_p}([\mathcal{G}_n]) = i_* d_n \quad \mathrm{in} \quad A_* (Y_p)_{\mathbb{Q}}
\end{equation}
for some $d_n \in A_* (b_p)_{\mathbb{Q}}$.  Let $f:Y_p \to \mathrm{Spec}(k)$
be the structure morphism.  Then applying $f_*$ to both sides
of (\ref{eq:pushit}) and using the Riemann-Roch result of 
 \cite[Thm. 18.2]{Ful} leads to
 \begin{equation}
\label{eq:easier}
\mathrm{dim}_k (\HH^0(Y_p,(\pi_p)_*\pi_p^* \mathcal{G}_n)  =
\# G_p \cdot \mathrm{dim}_k (\HH^0(Y_p, \mathcal{G}_n)) + f_* i_* d_n\quad \mathrm{in}\quad 
A_* (\mathrm{Spec}(k))_{\mathbb{Q}} = \mathbb{Q}
\end{equation}
for $n>>0$.

Now write $n = a + rm$ as before.   Then
$\mathcal{G}_n = \mathcal{G}_{a+rm} = \mathcal{G}_a\otimes (\zeta^* \mathcal{L}_Y)^r$.
Since $\tau_{Y_p}:G_0(Y_p) \to A_* (Y_p)_{\mathbb{Q}}$ is compatible
with respect to cap products, we find from part (2) of 
\cite[Thm. 18.2]{Ful} that
$$i_* d_{a+rm} = i_*( \mathrm{ch}_b(\mathcal{L}_b^r) \cdot d_a ) = i_* \left (\sum_{j = 0}^{c_p} \frac{c_1(\mathcal{L}_b)^j r^j }{j!} \cdot d_a\right )$$
where $\mathcal{L}_b = i^* \zeta^* \mathcal{L}_Y$, 
$\mathrm{ch}_b$ is the Chern character associated to vector bundles on
$b_p$, $c_1(\mathcal{L}_b)$ is the first Chern class of $\mathcal{L}_b$ and $\mathrm{dim}(b_p) =  c_p$.  Since $f_*: A_* (Y_p)_{\mathbb{Q}} \to A_* (\mathrm{Spec}(k))_{\mathbb{Q}}$ is linear, it follows that the function
$r \to f_* i_* d_{a+rm}$ is a polynomial of degree bounded by $c_p$ in $r$.  Thus  (\ref{eq:easier}) 
implies (\ref{eq:zapzap}), which completes the proof of (\ref{eq:littleproj})
and also shows that $P(S(X,\mathcal{F},\mathcal{L}))$ has polynomial growth.
Finally, the equality $\mathrm{max}\{c_p,d(P(S(X,\mathcal{F},\mathcal{L})))\}
=d$ follows from the fact that $S(X,\mathcal{F},\mathcal{L})
=P(S(X,\mathcal{F},\mathcal{L}))\oplus P'(S(X,\mathcal{F},\mathcal{L}))$.


\subsection{Proof of part (iii) of Theorem \ref{thm:meow}}
\label{ss:meow3}

We assume the notations of the previous two subsections.    
In view of Lemma \ref{lem:nolabel},
to complete
the proof of Theorem \ref{thm:meow} it will suffice to show the following.  
Let $g \in G$ be a $p$-regular element which is not the identity element.
Define $B$ to be the ramification locus of the cover $\pi:X\to Y$, and let
$c=\mathrm{dim}(\mathrm{supp}(\mathcal{F})\cap B)$.
Suppose $0 \le a < \#G$ and that $n = a + r \,\#G$ for some $0 \le r \in \mathbb{Z}$.  It will
suffice to show that there is a polynomial 
$\rho_{g,a}(t)  \in F[t]$ of degree less than or equal to $c$
in $t$ such that 
\begin{equation}
\label{eq:phiequal}
\Phi(\HH^0(X,\mathcal{F} \otimes \mathcal{L}^{a+r \,\#G}))(g) = \rho_{g,a}(r)\quad \mathrm{for}\quad r >> 0.
\end{equation}
By replacing $X$ by $\mathrm{supp}(\mathcal{F})$, we can reduce to the case
in which $c = \mathrm{dim}(B)$.

We will construct $\rho_{g,a}(t)$ using a result of Baum, Fulton and Quart \cite{BFQ}.
Let $G_0(G,X)$  be the Grothendieck group of 
coherent sheaves on $X$ having a compatible
action of $G$.   Let $H$ be the subgroup of $G$ generated by $g$. The subscheme $X^H$ of $X$ fixed by $H$ is closed and hence projective since $X$ is
projective. Since $k$ is algebraically closed and $g$ is a $p$-regular 
element, $kH$ is the direct sum of idempotent subspaces indexed by 
the character 
group $\hat {H}=\mathrm{Hom}(H,k^*)$ of $H$.  This gives rise to an  isomorphism 
$$G_0(H,X^H) = G_0(X^H) \otimes_{\mathbb{Z}} G_0(kH)$$
where $G_0(kH)$ is isomorphic  to $\mathbb{Z}\hat{H}$.  Let $\Phi_g:G_0(kH) \to \Lambda$ be the ring homomorphism which sends the class $[M]$ of a $kH$-module
$M$ to $\Phi(M)(g)$.   We then have a homomorphism
\begin{equation}
\label{eq:G0map}
1 \otimes \Phi_g: G_0(H,X^H) = G_0(X^H) \otimes_{\mathbb{Z}} G_0(kH) \to G_0(X^H) \otimes_{\mathbb{Z}} \Lambda.
\end{equation}

Baum, Fulton and Quart show there is a commutative diagram
\begin{equation}
\label{eq:crazytwo}
\xymatrix @-.8pc {
 G_0(G,X) \ar[rrrr]^(.5){L_*} \ar@<2ex>[d]_{\chi^{\mathrm{naive}}} &&&&
G_0(X^H) \otimes_{\mathbb{Z}} \Lambda \ar[d]^{\chi^{\mathrm{dim}} \otimes 1}  \\
 G_0(kG) \ar[rrrr]^(.5){\Phi_g} &&&& \mathbb{Z} \otimes_{\mathbb{Z}}\Lambda  = \Lambda 
}
\end{equation}
in which $\chi^{\mathrm{naive}}:G_0(G,X)\to G_0(kG)$ (resp. 
$\chi^{\mathrm{dim}}:G_0(X^H) \to \mathbb{Z}$) is the natural
equivariant (resp. non-equivariant) Euler characteristic homomorphism, 
and the homomorphism $L_*$ has the following property.
Let $K_0(H,X^H)$
be the Grothendieck group of all locally free coherent sheaves on $X^H$
with a compatible action of $H$.  Then $K_0(H,X^H)$ is a ring via the
tensor product
over $\mathcal{O}_{X^H}$, and $G_0(H,X^H)$ is a module
for this ring. For all locally
free $G$-sheaves $\mathcal{E}$ on $X$ we have 
\begin{equation}
\label{eq:cap}
L_*([\mathcal{E} \otimes \mathcal{F}]) = [\mathrm{res}(\mathcal{E})] \cdot 
L_* ([\mathcal{F}])
\end{equation}
where $[\mathrm{res}(\mathcal{E})]$ is the class in $K_0(H,X^H)$ determined
by the restriction $\mathrm{res}(\mathcal{E})$ of $\mathcal{E}$
to $X^H$.  

We now apply diagram (\ref{eq:crazytwo}) to the class $[\mathcal{E} \otimes \mathcal{F}] \in G_0(G,X)$
where $\mathcal{E} = \mathcal{L}^n$,  $n = a + r \,\#G$
and $0 \le r \in \mathbb{Z}$.  The restriction $\mathrm{res}(\mathcal{L})$ of $\mathcal{L}$ to
$X^H$ is a line bundle on $X^H$, with a compatible action of $H = \langle g \rangle$.
Since $\mathrm{res}(\mathcal{L})$ is the direct sum of its images under
the idempotents of $kH$, we conclude that for each stalk $\mathrm{res}(\mathcal{L})_P$ of 
$\mathrm{res}(\mathcal{L})$ at a point $P \in X^H$, all but one of these idempotents
annihilates $\mathrm{res}(\mathcal{L})_P$, and the action of $H$ on $\mathrm{res}(\mathcal{L})_P$
is via a single character in $\hat{H}$.  
Hence $\mathrm{res}(\mathcal{L}^{\# G})$ is a line bundle on $X^H$ with trivial $H$-action.
So in $K_0(H,X^H) = K_0(X^H) \otimes_{\mathbb{Z}} G_0(kH)$ we have
\begin{equation}
\label{eq:equaler}
[\mathrm{res}(\mathcal{E})] = [\mathrm{res}(\mathcal{L}^{a + r \,\# G})] = 
(\epsilon^r \otimes 1) \cdot [\mathrm{res}(\mathcal{L}^{a })]
\end{equation}
where $\epsilon = [\mathrm{res}(\mathcal{L}^{\# G})] \in K_0(X^H)$
and $(\epsilon^r \otimes 1) \in K_0(X^H) \otimes_{\mathbb{Z}} \Lambda$. 
Since $\mathcal{L}$ is ample, we have $\HH^i(X,\mathcal{F} \otimes \mathcal{L}^{a + r \,\# G}) = 0$
for all $i > 0$ if $r$ is sufficiently large.  Thus $\chi^{\mathrm{naive}}([\mathcal{E} \otimes\mathcal{F}]) 
= [\HH^0(X,\mathcal{E} \otimes \mathcal{F})]$ in $G_0(kG)$, so
$$\Phi_g(\chi^{\mathrm{naive}}([\mathcal{E} \otimes \mathcal{F}])) = 
\Phi(\HH^0(X,\mathcal{F} \otimes \mathcal{L}^{a + r \,\# G}))(g)$$
for $t >> 0$.  Hence going the other way around the diagram in (\ref{eq:crazytwo})
and using (\ref{eq:equaler}) gives
\begin{equation}
\label{eq:otherway}
\Phi(\HH^0(X,\mathcal{F} \otimes \mathcal{L}^{a + r \,\# G}))(g) 
= (\chi^{\mathrm{dim}}\otimes 1) \left ( (\epsilon^r \otimes 1) \cdot  
[\mathrm{res}(\mathcal{L}^{a })] \cdot L_* ([\mathcal{F}]) \right )
\end{equation}
for $r >>0$.
Since $\epsilon$ is the class in $K_0(X^H)$ of the ample line bundle $\mathrm{res}(\mathcal{L}^{\# G})$ on $X^H$, and $X^H$ is a projective scheme over $k$, the existence of Hilbert polynomials implies that the function
 of $r$ on the right hand side of (\ref{eq:otherway}) is for large $r$ given by a polynomial in $r$ (with coefficients in $F$) of degree bounded by $\mathrm{dim}(X^H)$.  Since $H$ is not the trivial subgroup of $G$,
 $X^H$ is contained in the ramification locus $B$.  
Thus $\mathrm{dim}(X^H) \le  \mathrm{dim}(B) = \mathrm{dim}(B \cap \mathrm{supp}(\mathcal{F})) = c$.  This shows
 that there is a polynomial $\rho_{a,g}(t) \in F[t]$ of the sort required  in (\ref{eq:phiequal}),
 which completes the proof.
 

\subsection{Proof of Corollary \ref{cor:ohyeah}}
\label{ss:meowcor}

We assume the notations of the previous subsections. In particular, $\pi:X \to Y$, 
$\pi_p:X \to Y_p = X/G_p$ and $\zeta:Y_p \to Y$ are the natural  
quotient morphisms, and $B_p$ denotes the ramification locus of $\pi_p$. 
As before, we may reduce
to the case $X=\mathrm{supp}(\mathcal{F})$.
If $\mathrm{char}(k)=p=0$ or $B_p=\mathrm{supp}(\mathcal{F})\cap B_p=\emptyset$,
Corollary \ref{cor:ohyeah} follows from Theorem \ref{thm:specialtame}.

We now consider the remaining case of Corollary \ref{cor:ohyeah}. Hence 
we suppose $\mathrm{char}(k)=p>0$ and 
$\mathrm{dim}(B_p)=\mathrm{dim}(\mathrm{supp}(\mathcal{F})\cap B_p)=0$.

\begin{sublem}
\label{lem:bigdegree}
To show Corollary $\ref{cor:ohyeah}$, it will suffice to show that when $k$
is algebraically closed,  there is an integer $m > 0$ such that
there is an injection for all sufficiently large $n$
\begin{equation}
\label{eq:into1}
P'(\HH^0(X,\mathcal{F}_n)) \to 
P'(\HH^0(X,\mathcal{F}_{n+m}))
\end{equation}
when $\mathcal{F}_n=\mathcal{F}\otimes\mathcal{L}^n$. 
\end{sublem}

\begin{proof}
By Lemma \ref{lem:biggerfield}, we can reduce
to the case in which $k$ is algebraically closed.  Theorem \ref{thm:meow} shows 
$\mathrm{dim}_k(P'(H^0(X,\mathcal{F}_n)))$ 
is bounded independently of $n$.  Hence given an $m$ as in the statement of the 
Lemma, it follows that for sufficiently large $n$,
the isomorphism class of $P'(H^0(X,\mathcal{F}_n))$
depends only on $n\mod m$.  Hence 
$P'(S(X,\mathcal{F},\mathcal{L}))$ has a polynomial 
description of degree $0$.  By Theorem \ref{thm:thebigtwo}, 
$S(X,\mathcal{F},\mathcal{L})$ has a $G_0(kG)$-polynomial description   of degree bounded by $d$.  Since 
$S(X,\mathcal{F},\mathcal{L}) = P(S(X,\mathcal{F},\mathcal{L})) \oplus
P'(S(X,\mathcal{F},\mathcal{L}))$,
it follows that   $P(S(X,\mathcal{F},\mathcal{L}))$
has a $G_0(kG)$-polynomial description of degree bounded by $d$.  Because 
the terms of  $P(S(X,\mathcal{F},\mathcal{L}))$ are projective
and finitely generated, and $K_0(kG) \to G_0(kG)$ is injective, we conclude that both
$P(S(X,\mathcal{F},\mathcal{L}))$ and
$S(X,\mathcal{F},\mathcal{L})$ have polynomial descriptions of
degree bounded by $d$.  The degrees are exactly $d$ because the
ordinary Hilbert polynomial of $S(X,\mathcal{F},\mathcal{L})$ has degree
$d = \mathrm{dim}(\mathrm{supp}(\mathcal{F}))$.  This implies 
$S(X,\mathcal{F},\mathcal{L})$ is indecomposably finite.
\end{proof}

The trace operator $\mathrm{Tr}_{G_p}$ gives an exact sequence
\begin{equation}
\label{eq:tracit}
0 \to \mathrm{Tr}_{G_p} (\pi_{p})_* (\mathcal{F}_n) \to  
(\pi_{p})_* (\mathcal{F}_n) \to \mathcal{C}_n \to 0
\end{equation}
of sheaves on $Y_p$, where $\mathcal{C}_n$ is coherent.

\begin{sublem}
\label{lem:prodinj}
To produce an injection $(\ref{eq:into1})$, it will suffice to show there is an integer $m>0$ and a very ample line bundle $\mathcal{L}_Y$ on $Y = X/G$ with the following properties:
\begin{enumerate}
\item[i.] There is a $G$-equivariant
$O_X$-module isomorphism $\mathcal{L}^m = \pi^* \mathcal{L}_Y$.  
\item[ii.] There is an injection  $\mathcal{O}_{Y_p} \to \zeta^* \mathcal{L}_Y   $ whose tensor product with  $(\ref{eq:tracit})$ gives a commutative diagram
\begin{equation}
\label{eq:arraydiag}
\xymatrix @-1pc {
0\ar[r]& \mathrm{Tr}_{G_p} (\pi_{p})_* (\mathcal{F}_n)\ar[r]\ar[d]&
(\pi_{p})_* (\mathcal{F}_n)\ar[r]\ar[d]&\mathcal{C}_n \ar[r]\ar[d]&0\\
0\ar[r]& \mathrm{Tr}_{G_p} (\pi_{p})_* (\mathcal{F}_n)\otimes \zeta^* \mathcal{L}_Y\ar[r]&
(\pi_{p})_* (\mathcal{F}_n)\otimes \zeta^* \mathcal{L}_Y\ar[r]&\mathcal{C}_n \otimes \zeta^* \mathcal{L}_Y\ar[r]&0
}
\end{equation}
in which the vertical morphisms are injective.
\end{enumerate}
\end{sublem}

\begin{proof} 
The projection formula together with the fact that  $\pi^* \mathcal{L}_Y \cong 
\mathcal{L}^m$ shows
that the middle vertical morphism in (\ref{eq:arraydiag}) is a morphism 
\begin{equation}
\label{eq:itsinj}
(\pi_{p})_* (\mathcal{F}_n) \to 
(\pi_{p})_* (\mathcal{F}_n)\otimes \zeta^* \mathcal{L}_Y \cong
(\pi_{p})_*(\mathcal{F}_{n+m}).
\end{equation}
Taking cohomology gives a morphism
$$\tau_n:\HH^0(X,\mathcal{F}_n) = 
\HH^0(Y_p,(\pi_{p})_*(\mathcal{F}_n)) \to 
\HH^0(Y_p,(\pi_{p})_*(\mathcal{F}_{n+m})) = \HH^0(X,
\mathcal{F}_{n+m}).$$
Since (\ref{eq:itsinj}) is injective by hypothesis, $\tau_n$ is injective. 
So to prove that $\tau_n$ gives an injection of the form in
(\ref{eq:into1}), it will suffice to show that 
\begin{equation}
\label{eq:duhduhduh}
\tau_n(P'(\HH^0(X,\mathcal{F}_n))) \cap 
P(\HH^0(X,\mathcal{F}_{n+m})) = \{0\}\quad \mbox{ for }n>>0
\end{equation}
for any choice of maximal projective summand $P(\HH^0(X,\mathcal{F}_{n+m}))$
of $\HH^0(X,\mathcal{F}_{n+m})$.  

The socle
of $P(\HH^0(X,\mathcal{F}_{n+m}))$ as a $kG_p$-module
is $\mathrm{Tr}_{G_p} P(\HH^0(X,\mathcal{F}_{n+m}))$.
Note that
$$P'(\HH^0(X,\mathcal{F}_n)) \cap 
\mathrm{Tr}_{G_p} \HH^0(X,\mathcal{F}_{n}) = \{0\}$$
since  $\mathrm{Tr}_{G_p}$ annihilates $P'(\HH^0(X,\mathcal{F}_n))$. So 
to prove (\ref{eq:duhduhduh}) it will suffice to show
 \begin{equation}
 \label{eq:duh7}
 \tau_n(\HH^0(X,\mathcal{F}_n)) \cap 
 \mathrm{Tr}_{G_p} \HH^0(X,\mathcal{F}_{n+m}) =
 \tau_n(\mathrm{Tr}_{G_p} \HH^0(X,\mathcal{F}_{n}))\quad \mbox{ for }n>>0.
 \end{equation}
 
 It was proved in (\ref{eq:asymptot}) that
 $$\mathrm{Tr}_{G_p} \HH^0(X,\mathcal{F}_{n}) = 
 \HH^0(Y_p,\mathrm{Tr}_{G_p} (\pi_p)_* \mathcal{F}_{n})\quad\mbox{ for }n>>0.$$
Hence, for $n>>0$, the cohomology over $Y_p$ of (\ref{eq:arraydiag}) gives a commutative diagram
\begin{equation}
\label{eq:arraydiag2}
\xymatrix @-1pc {
0\ar[r]& \mathrm{Tr}_{G_p} \HH^0(X,\mathcal{F}_{n})\ar[r]\ar[d]&
\HH^0(X,\mathcal{F}_n)\ar[r]\ar[d]&\HH^0(Y_p,\mathcal{C}_n) \ar[d]\\
0\ar[r]&  \mathrm{Tr}_{G_p} \HH^0(X,\mathcal{F}_{n+m})\ar[r]&
\HH^0(X,\mathcal{F}_{n+m})\ar[r]&\HH^0(Y_p,\mathcal{C}_n \otimes \zeta^* \mathcal{L}_Y)
}
\end{equation}
in which the vertical morphisms are induced by $\tau_n$.  Therefore (\ref{eq:duh7})
follows from the fact that the right vertical morphism in (\ref{eq:arraydiag2}) is injective
due to the fact that the right vertical morphism in (\ref{eq:arraydiag}) has been
assumed to be injective.
\end{proof}

To construct a line bundle $\mathcal{L}_Y$ with the properties in Lemma \ref{lem:prodinj},
we will use the following Lemma.

\begin{sublem}
\label{lem:support}
Suppose $Z$ is a projective scheme over an algebraically closed field $k$
and that $\mathcal{T}$ is a coherent sheaf on $Z$.   There is a finite filtration of $\mathcal{T}$ by
coherent sheaves such that the successive quotients are isomorphic
to line bundles on a finite set $\{Z_i\}_i$ of integral closed subschemes of $Z$.  
Let $E$ be an effective
Cartier divisor on $Z$.   Then the natural morphism $\mathcal{T} \to \mathcal{T} \otimes \mathcal{O}_Z(E)$ is injective provided that the support of $E$ does not contain
the generic point of any of the $Z_i$.
\end{sublem}

\begin{proof}  The first statement follows from \cite[Prop. I.7.4]{hart}.  Since
tensoring with $\mathcal{O}_Z(E)$ over $\mathcal{O}_Z$ is exact, to prove
the second statement we can reduce to the case in which $\mathcal{T}$ is
a line bundle on one of the $Z_i$.  We have an exact sequence of sheaves
$$\mathrm{Tor}^1_{\mathcal{O}_Z}(\mathcal{T},\mathcal{O}_Z(E)/\mathcal{O}_Z) \to \mathcal{T} \to \mathcal{T} \otimes \mathcal{O}_Z(E).$$
This shows that the kernel $\mathcal{K}$ of $\mathcal{T} \to \mathcal{T} \otimes \mathcal{O}_Z(E)$
is annihilated by the annihilator ideal of $\mathcal{O}_Z(E)/\mathcal{O}_Z$.
Hence the support of $\mathcal{K}$ is contained in the support of $E$.
Since $\mathcal{K}$ is a subsheaf of the line bundle  $\mathcal{T}$ 
on the integral scheme $Z_i$,  the support of $\mathcal{K}$ is all of $Z_i$
if $\mathcal{K}$ is not $0$.  However, the support of $E$
does not contain the generic point of $Z_i$ by assumption, so $\mathcal{K}$
must be $0$.
\end{proof}

\medskip

\noindent
\textit{Completion of the proof of Corollary $\ref{cor:ohyeah}$.}
By Lemma \ref{lem:pullback}, there is a very ample line bundle $\mathcal{L}_Y$
on $Y$ such that $\pi^* \mathcal{L}_Y$ is $G$-isomorphic
to $\mathcal{L}^m$ for some integer $m>0$.   We
we can use Bertini's theorem over the algebraically closed field $k$ to 
move $\mathcal{L}_Y$ so that it is defined by an effective Cartier divisor $D$ on $Y$ with the following property.
Suppose $0 \le a < m$ and that $\mathcal{T}$ is one of the sheaves
$\mathrm{Tr}_{G_p}(\pi_{p})_*(\mathcal{F}_{a})$,
$(\pi_{p})_* ( \mathcal{F}_{a})$ or $\mathcal{C}_a$
on $Y_p$.  If $Z_i$ is an integral subscheme of $Y_p$ arising
from a filtration of $\mathcal{T}$ of the kind described in Lemma \ref{lem:support},
then $D$ does not contain the image in $Y$ of the generic point of $Z_i$.
Observe that since there are only finitely many $a$ involved, these are finitely
many conditions on $D$.  Because $\mathcal{L}^m \cong \pi^* \mathcal{L}_Y$,
replacing $a$ in the above condition by $a + rm$ for some $r \ge 0$ only changes
each $\mathcal{T}$ by tensoring it with a power of $\zeta^*\mathcal{L}_Y$.  
This only changes
the line bundles on the $Z_i$ which arise in Lemma \ref{lem:support},
so we conclude that $D$ has the properties above when $a$ is any positive
integer.  We now conclude from Lemma \ref{lem:support} that $\mathcal{L}_Y$ has
the properties needed in Lemma \ref{lem:prodinj}, which completes
the proof. \hfill $\BBox$


\section{Indecomposable finiteness for smooth projective varieties over fields}
\label{s:projvarfields}
\setcounter{equation}{0}

The object of this section is to prove Theorem \ref{thm:main}.  
Throughout this section we will assume the notations of \S \ref{s:hom}. Moreover, we will assume that the ring $R$ of \S \ref{s:hom} is a field $k$ of positive characteristic $p$, and that the action of $G$ on $X$ is faithful and generically free
over $k$.

\subsection{Line bundles and schemes of dimension $0$ or $1$}
\label{ss:prelimvar}

By Lemma \ref{lem:biggerfield}, we can assume the following hypothesis:

\begin{subhypo}
\label{hypo:whichscheme} The scheme $X$ is a smooth projective variety 
of dimension $d$ over
an algebraically closed field $k$  of characteristic $p > 0$, and 
the action of $G$ on $X$ is faithful and generically free over $k$.
Let $B\subset X$ be the ramification locus
of the quotient morphism $\pi:X \to Y=X/G$, and let $b$ the the image of $B$ in $Y$.
The $G$-equivariant line bundle $\mathcal{L}$ on $X$ is ample.
\end{subhypo}

Note that $Y=X/G$ is a Noetherian normal variety.

\begin{subnot}
\label{notation}
In this section, we will write $S(X,\mathcal{L})$ for $S(X,\mathcal{O}_X,\mathcal{L})$;
and $\HH^i(\mathcal{F})$ for $\HH^i(X,\mathcal{F})$ if 
$\mathcal{F}$ is a coherent sheaf on $X$.
\end{subnot}

We need the following refinement of Lemma \ref{lem:pullback}.

\begin{sublem}
\label{lem:linebundle}
 Let $A$ be
a finite set of closed points of $X$.
Let $m_1$ be a positive integer such that $\mathcal{L}^{m_1}$ is very ample. 
Let $m_2$ be the maximal divisor of $\# G$ which is prime to $p$.
Then there exist a $G$-stable effective very ample Weil divisor $T$ on $X$ and effective Weil divisors $D_0,\ldots,D_{d-1}$ on $Y$ such that when
$D_j'=\pi^*(D_j)$ the following is true:
\begin{enumerate}
\item[i.] There is a $G$-equivariant $\mathcal{O}_X$-module
isomorphism $\psi:\mathcal{L}^{m_1  m_2^2}\to\mathcal{O}_X(T)$;
\item[ii.] for all $j$, $A\cap D_j'=\emptyset$;
\item[iii.] there is a $G$-equivariant $\mathcal{O}_X$-module isomorphism 
$$\lambda_j:\quad  \mathcal{L}^{m_1 m_2^2  (\#G)}\to \mathcal{O}_X(D_j');$$
\item[iv.] the divisors $D_0,\ldots,D_{d-1}$ are all linearly equivalent on $Y$;
\item[v.] for each subset $J\subset \{0,\ldots,d-1\}$ and each point 
$x \in \cap_{j \in J} D'_j$, the local equations of the $D'_j$ in 
$\mathcal{O}_{X,x}$ for $j \in J$ form
a regular sequence; 
\item[vi.] the intersection $D_0'\cap \cdots \cap D_{d-1}'\cap B=\emptyset$.
\end{enumerate}
\end{sublem}

\begin{proof}
Let $m_1$ be a positive integer such that $\mathcal{L}^{m_1}$ is very ample.
Then $\HH^0(\mathcal{L}^{m_1})$ contains a non-trivial finite $G$-stable subgroup
of order a power of $p$.  Hence there is a non-zero element   $f \in \HH^0(\mathcal{L}^{m_1})$ whose orbit $\omega(f)$ under the action of $G$ has order 
$m_3$ dividing $m_2$.    The element
$h = (\otimes_{\tau \in \omega(f)}\, \tau)^{m_2/m_3}$ defines a global section of 
$\mathcal{L}^{m_1 m_2}$ which is non-zero since each $\tau \in \omega(f)$ is non-zero.
Consider the $G$-equivariant embedding $X\hookrightarrow \PP_k^{N'}$
associated to $\HH^0(\mathcal{L}^{m_1 m_2})$.  The section $h$ defines a $G$-stable hyperplane $T'$ in $\PP_k^{N'}$, and we let
$T=T'\cap X$. Note that since $h\neq 0$ on $X$, we ensure that $T'$ does not
contain $X$. We obtain an $\mathcal{O}_X$-module isomorphism $\psi_1:\mathcal{L}^{m_1 m_2}\to
\mathcal{O}_X(T)$. There is a $G$-action on $\mathrm{Hom}_{\mathcal{O}_X}(\mathcal{L}^{m_1 m_2}, 
\mathcal{O}_X(T))$ given by $(g\,\rho)(t)=g\, \rho(g^{-1}t)$ for all $\rho\in 
\mathrm{Hom}_{\mathcal{O}_X}(\mathcal{L}^{m_1 m_2},\mathcal{O}_X(T))$ and
all $g\in G$. Hence for all $g\in G$, $(g\,\psi_1): \mathcal{L}^{m_1 m_2} \to 
\mathcal{O}_X(T)$ is an $\mathcal{O}_X$-module isomorphism. We obtain a one-cocycle
from $G$ to $\mathrm{Aut}_{\mathcal{O}_X}(\mathcal{O}_X(T))$ via $g\mapsto \psi_1\circ
(g\,\psi_1)^{-1}$. Since $k$ is algebraically closed and $X$ is a projective variety over $k$,
$\mathrm{Aut}_{\mathcal{O}_X}(\mathcal{O}_X(T))={k}^*$. Thus the one-cocycle above defines a class 
in $\HH^1(G,{k}^*)$, where $G$ acts trivially on $k^*$.  This cohomology class is annihilated
by $m_2$. So  on replacing $\mathcal{L}^{m_1 m_2}$ by $\mathcal{L}^{m_1 m_2^2}$, 
$h$ by $h^{m_2}$, and $T$ by the resulting divisor on $X$, we can assume the above class  in $\HH^1(G,k^*)$ is trivial.  This means that we can change
$\psi_1$ by a scalar in ${k}^*$ so as to obtain a $G$-equivariant $\mathcal{O}_X$-module
isomorphism $\psi:\mathcal{L}^{m_1 m_2^2}\to\mathcal{O}_X(T)$, which implies (i).
We let $X\hookrightarrow \PP_k^{N}$ be the $G$-equivariant embedding 
associated to $\HH^0(\mathcal{L}^{m_1 m_2^2})$.

Suppose $r_0,\ldots,r_{d-1}$ are linear forms on $\PP_k^N$.  
For each $1 \le l \le d$ and each $l$-tuple 
 $(j_1,\ldots,j_l)$ of distinct integers in $\{0,\ldots,d-1\}$ and each $l$-tuple
 $(g_1,\ldots,g_l)$ of elements of $G$, let $L(g_1 r_{j_1},\ldots,g_l r_{j_l})$
 be the common zero locus of $g_1 r_{j_1},\ldots,g_l r_{j_l}$ in $\PP_k^N$.
 Define
 $$\Theta(g_1 r_{j_1},\ldots,g_l r_{j_l}) = L(g_1 r_{j_1},\ldots,g_l r_{j_l}) \cap X.$$
 For $0 \le j \le d-1$ let
 $$D'_j = \sum_{g \in G} \Theta(gr_j).$$
 Since $k$ is an algebraically closed field, we can use Bertini's Theorem
 to find $r_0,\ldots,r_{d-1}$ so the following is true.
 \begin{enumerate}
 \item[1.] $\Theta(r_j) \cap \bigcup_{g \in G} gA = \emptyset$ for all $j$.
 \item[2.] $\Theta(r_j)$ does not contain an irreducible component of $T$ for any $j$.
 \item[3.] Fix $j \in \{0,\ldots,d-1\}$.  Suppose 
 $(j_1,\ldots,j_l)$ is an $l$-tuple of distinct integers in $\{0,\ldots,j-1\}$.  Then 
 $\Theta(r_j,g_1 r_{j_1},\ldots,g_l r_{j_l}) = \Theta(r_j) \cap \Theta(g_1 r_{j_1},\ldots,g_l r_{j_l})$ has codimension
 $l+1$ in $X$ for each $l$-tuple $(g_1,\ldots,g_l)$ of elements of $G$.
 Moreover, $\Theta(r_j, g_1 r_{j_1},\ldots,g_l r_{j_l})\cap B$ has codimension at least $l+2$ in $X$.
 \end{enumerate}
One applies Bertini's Theorem to show that once $r_0,\cdots,r_{j-1}$
 have been chosen, each of the conditions in (1) - (3) is satisfied by a non-empty open dense set of
 linear forms $r_j$ in the space of all linear forms.

We now consider the rational function
$\displaystyle \frac{\mathrm{Norm}(r_j)}{h^{\#G}} = \frac{\prod_{g \in G} gr_j}{h^{\# G}}$ on 
$\PP_k^N$, which restricts to 
a rational function $l_j$ on $X$. Then $l_j$ is a $G$-invariant rational function 
in $k(Y)=k(X)^G$ such that $\mathrm{div}_{Y}(l_j)$ pulls back to 
$$\mathrm{div}_{X}(l_j) = D'_j - \# G \cdot T.$$ 
Since $Y$ is a noetherian normal variety, it is regular in codimension 1,
and we define $D_j$ to be the Weil divisor associated to the effective part of 
$\mathrm{div}_{Y}(l_j)$.  Since no irreducible component of the $G$-invariant
divisor $T$ occurs in $D'_j$, $\pi^*(D_j)  = D'_j$ is the 
zero locus of $\mathrm{Norm}(r_j)$ on $X$.  Condition (3) implies that for
each subset  $J\subset \{0,\ldots,d-1\}$, the codimension of $\bigcap_{j \in J} D'_j$
in $X$ is equal to $\# J$.  By \cite[Thm. 31, p. 108]{Matsumura}, this implies condition (v) in
the statement of Lemma \ref{lem:linebundle}.   Defining $\lambda_j=\mathrm{mult}(l_j^{-1})\circ\psi^{\otimes \#G}$ leads to (iii).  Conditions (ii) and (iv) are clear.
For (vi) we use (3), which completes the proof.
\end{proof}

We will make use of the following result of Nakajima to split various exact sequences
of cohomology groups.

\begin{sublem}
\label{lem:project}
Suppose that $Z$ is a scheme of dimension $0$ over $k$ on which the finite group
$G$ acts freely over $k$.  Let $\mathcal{F}$ be a coherent $G$-sheaf on $Z$.  Then
$\HH^0(Z,\mathcal{F})$ is a free $kG$-module.
\end{sublem}

\begin{proof} By definition, the quotient morphism $\pi:Z \to Z/G$ is a $G$-torsor
when we view $G$ as a finite constant group scheme.  By descent, $\mathcal{F}$ is isomorphic
to $\pi^*(\mathcal{H})$ for some coherent sheaf $\mathcal{H}$ on $Z/G$.  By \cite[Theorem 1]{N},
there is a bounded complex $C^\bullet $ of finitely generated free 
$kG$-modules
such that $\HH^i(Z,\pi^*(\mathcal{H})) = \HH^i(Z,\mathcal{F})$ is isomorphic to $\HH^i(C^\bullet)$ for all $i$.
Since $Z$ has dimension $0$, these groups  vanish for $i \ne 0$.  Since
finitely generated free $kG$-modules are injective, we can split $C^\bullet$
completely, which implies Lemma \ref{lem:project}.
\end{proof}

We end this subsection by observing that part (i) of Theorem \ref{thm:main} concerning
curves $X$ follows from Corollary \ref{cor:ohyeah}.

\subsection{Smooth projective acyclic surfaces}
\label{ss:surface}

Let $X$ be a smooth projective acyclic surface such that
$X$,  $\mathcal{L}$, $Y$, $\pi$ and $B$
are as in Hypothesis $\ref{hypo:whichscheme}$.
(Recall that $X$ is acyclic if $\HH^i(\mathcal{O}_X)=0$ for all
$i\geq 1$.) Using Lemma \ref{lem:linebundle} we fix the following
notation.

\begin{subhypo}
\label{hypo:Lamp}The integer $m_1$ is large enough so that 
$\mathcal{L}^{m_1}$
is very ample and $\HH^i(X,\mathcal{L}^{m_1 s}) = 0$ if $i, s > 0$. 
Let  $m=\#G$
and let $m_2$ be the maximal divisor of $\# G$ which is prime to $p$.
Set $\mathcal{T} = \mathcal{L}^{m_1 m_2^2}$.  Let 
$A$ be a finite set of closed
points of $X$ to be specified later.  By Lemma \ref{lem:linebundle}, there are effective divisors $D$ and $E$
on the normal surface $Y = X/G$ with the following properties.  There is a function $f \in k(Y)$
such that $\mathrm{div}_Y(f)=E-D$.  Define $D' = \pi^*(D)$ and $E' = \pi^*(E)$.
There are $G$-equivariant $\mathcal{O}_X$-module
isomorphisms $\lambda: \mathcal{T}^m\to \mathcal{O}_X(D') =_{\mathrm{def}} \mathcal{T}_{D'}$,
and $\psi:\mathcal{T} \to \mathcal{O}_X(T)$ where $T$ is as in Lemma \ref{lem:linebundle}. 
At each point  $x \in X$ where $E'$ and $D'$ intersect, the local
equations of these divisors in $\mathcal{O}_{X,x}$ form a regular sequence.  Finally,
$\emptyset = E'\cap D'\cap B = E' \cap A = D' \cap A$. 
\end{subhypo}

\begin{sublem}
\label{lem:diagram}
Let $j$ be an integer with $0\leq j< m =\#G$. Then for all positive 
integers $t$, there exists a commutative diagram of sheaves of 
$G$-$\mathcal{O}_X$-modules
\begin{equation}
\label{eq:bigdiagram}
\xymatrix @-1pc {
&&0\ar[d]&0\ar[d]&\\
0\ar[r]&f\mathcal{T}_{D'}^{-(t+1)}\otimes \mathcal{T}^{(t+1)m+j}\ar[r]&
\mathcal{T}_{D'}^{-t}\otimes \mathcal{T}^{(t+1)m+j}\ar[r]\ar[d]&\mathcal{T}_{D'}^{-t}\otimes 
\mathcal{T}^{(t+1)m+j}\big|_{E'}\ar[r]\ar[d]&0\\
0\ar[r]&f\mathcal{T}_{D'}^{-1}\otimes \mathcal{T}^{(t+1)m+j}\ar[r]&
\mathcal{T}^{(t+1)m+j}\ar[r]\ar[d]&\mathcal{T}^{(t+1)m+j}\big|_{E'}\ar[r]\ar[d]&0\\
&&\mathcal{T}^{(t+1)m+j}\big|_{tD'}\ar[r]\ar[d]&\mathcal{T}^{(t+1)m+j}\big|_{tD'\cap E'}\ar[d]&\\
&&0&0&}
\end{equation}
which consists of $4$ short exact sequences.
\end{sublem}

\begin{proof}
The first row, resp.  second row, results from tensoring the short exact sequence
$$0\to \mathcal{O}_X(-E')\to \mathcal{O}_X \to \mathcal{O}_{E'}\to 0$$
with $\mathcal{T}_{D'}^{-t}\otimes \mathcal{T}^{(t+1)m+j}$, resp. with $\mathcal{T}^{(t+1)m+j}$, over 
$\mathcal{O}_X$. Note that the resulting two sequences stay exact, since 
$\mathcal{T}_{D'}^{-t}\otimes \mathcal{T}^{(t+1)m+j}$, resp. $\mathcal{T}^{(t+1)m+j}$, is a line bundle on $X$, hence a flat $\mathcal{O}_X$-module. We obtain the middle column by tensoring 
the short exact sequence
$$0\to \mathcal{O}_X(-tD')\to \mathcal{O}_X \to \mathcal{O}_{tD'}\to 0$$
with $\mathcal{T}^{(t+1)m+j}$ over $\mathcal{O}_X$. As before, the resulting sequence 
\begin{equation}
\label{eq:uhoh}
0\to \mathcal{T}_{D'}^{-t}\otimes \mathcal{T}^{(t+1)m+j}\to \mathcal{T}^{(t+1)m+j} \to 
\mathcal{T}^{(t+1)m+j}\big|_{tD'}\to 0
\end{equation}
stays exact. For the right column, we need to show that the short exact sequence (\ref{eq:uhoh}) 
stays exact when tensoring with $\mathcal{O}_{E'}$, i.e. that
\begin{equation}
\label{eq:uhoh2}
0\to\mathcal{T}_{D'}^{-t}\otimes \mathcal{T}^{(t+1)m+j}\big|_{E'}\arrow{\rho}{} 
\mathcal{T}^{(t+1)m+j}\big|_{E'} \to \mathcal{T}^{(t+1)m+j}\big|_{tD'\cap E'}\to 0
\end{equation}
is exact. It suffices to prove that $\rho$ is injective. 
Let $x\in E'\subset X$ and let $\delta_x$, resp. $\epsilon_x$,  be a local 
equation for $D'$, resp. $E'$,  at $x$. Let $\gamma_x$ be a local generator
for $\mathcal{T}$ at $x$, i.e. $\mathcal{T}_x=\gamma_x \cdot \mathcal{O}_{X,x}$. 
Thus, when restricting to the stalks at $x$, $\rho$ in (\ref{eq:uhoh2}) becomes 
$$\xymatrix @-1pc {
[\mathcal{T}_{D'}^{-t}\otimes \mathcal{T}^{(t+1)m+j}\big|_{E'}]_x\ar[r]^(.57){\rho_x} \ar@{=}[d]&
[\mathcal{T}^{(t+1)m+j}\big|_{E'}]_x\ar@{=}[d]\\
\frac{\delta_x^{t}\gamma_x^{\otimes((t+1)m+j)}\cdot \mathcal{O}_{X,x}}{\epsilon_x\cdot
\delta_x^{t}\gamma_x^{\otimes((t+1)m+j)}\cdot \mathcal{O}_{X,x}}&
\frac{\gamma_x^{\otimes((t+1)m+j)}\cdot \mathcal{O}_{X,x}}{\epsilon_x\cdot 
\gamma_x^{\otimes((t+1)m+j)}\cdot \mathcal{O}_{X,x}}
}$$
Hence $\rho_x$ is injective if and only if 
$$\delta_x^{t}\gamma_x^{\otimes((t+1)m+j)}\cdot \mathcal{O}_{X,x} \,\bigcap \,
\epsilon_x\cdot \gamma_x^{\otimes((t+1)m+j)}\cdot \mathcal{O}_{X,x}\subseteq 
\epsilon_x\cdot\delta_x^{t}\gamma_x^{\otimes((t+1)m+j)}\cdot \mathcal{O}_{X,x}.$$
Since $\gamma_x^{\otimes((t+1)m+j)}\cdot \mathcal{O}_{X,x}$ 
is a locally free $\mathcal{O}_{X,x}$-module and 
$\mathcal{O}_{X,x}$ is a unique factorization domain, to prove that
$\rho_x$ is injective, it suffices to prove that $\delta_x$ and $\epsilon_x$ 
are relatively prime elements of $\mathcal{O}_{X,x}$. 
This is clear, since at each point  $x \in E'\cap D'$,  the local equations $\delta_x$ and $\epsilon_x$ form a regular sequence in $\mathcal{O}_{X,x}$. 
Hence (\ref{eq:uhoh2}) is exact, which proves Lemma \ref{lem:diagram}.
\end{proof}

\begin{sublem}
\label{lem:cohomdiagram}
Let $j$ be an integer with $0\leq j< m= \#G$.  For all  positive integers $t$, we obtain from
the cohomology on $X$ of the terms of the diagram $(\ref{eq:bigdiagram})$ a commutative diagram of $kG$-modules
\begin{equation}
\label{eq:bigcohomdiagram}
\xymatrix @-1.2pc @R1pc {
&&0\ar[d]&0\ar[d]&\\
0\ar[r]&\HH^0(f\mathcal{T}_{D'}^{-(t+1)}\otimes \mathcal{T}^{(t+1)m+j})\ar[r]&
\HH^0(\mathcal{T}_{D'}^{-t}\otimes \mathcal{T}^{(t+1)m+j})\ar[r]\ar[d]&
\HH^0(\mathcal{T}_{D'}^{-t}\otimes \mathcal{T}^{(t+1)m+j}\big|_{E'})\ar[r]\ar[d]&0\\
0\ar[r]&\HH^0(f\mathcal{T}_{D'}^{-1}\otimes \mathcal{T}^{(t+1)m+j})\ar[r]&
\HH^0(\mathcal{T}^{(t+1)m+j})\ar[r]\ar[d]&\HH^0(\mathcal{T}^{(t+1)m+j}\big|_{E'})\ar[r]\ar[d]&0\\
&&\HH^0(\mathcal{T}^{(t+1)m+j}\big|_{tD'})\ar[r]^(.45){\tau}\ar[d]&\HH^0(\mathcal{T}^{(t+1)m+j}
\big|_{tD'\cap E'})\ar[d]&\\
&&0&0&}
\end{equation}
which consists of $4$ short exact sequences such that the homomorphism $\tau$ is surjective. \end{sublem}

\begin{proof}
By the Snake Lemma, it is enough to prove that the $4$ sequences are exact. We use the exactness of the corresponding $4$ sequences of $\mathcal{O}_X$-modules in (\ref{eq:bigdiagram}). By Hypothesis \ref{hypo:Lamp},  $\HH^1(\mathcal{T}^j)=0$  (when $j = 0$, recall that $X$ is acyclic).
Hence 
$$\HH^1(f\mathcal{T}_{D'}^{-(t+1)}\otimes \mathcal{T}^{(t+1)m+j})\cong
\HH^1(\mathcal{T}^{-(t+1)m}\otimes\mathcal{T}^{(t+1)m+j})
=\HH^1(\mathcal{T}^j)=0,$$
which gives the exactness of the first row of (\ref{eq:bigcohomdiagram}). 
Similarly,
$$\HH^1(f\mathcal{T}_{D'}^{-1}\otimes \mathcal{T}^{(t+1)m+j})\cong
\HH^1(\mathcal{T}^{tm+j})=0$$
and
$$\HH^1(\mathcal{T}_{D'}^{-t}\otimes \mathcal{T}^{(t+1)m+j})\cong
\HH^1(\mathcal{T}^{m+j})=0.$$
Hence the second row and the middle column of (\ref{eq:bigcohomdiagram}) are exact. To 
obtain that the right column of  (\ref{eq:bigcohomdiagram}) is exact, it suffices to show that 
$\HH^1(\mathcal{T}_{D'}^{-t}\otimes \mathcal{T}^{(t+1)m+j}\big|_{E'})=0$. 
By taking cohomology on $X$ of the short exact sequence
$$0\to f\mathcal{T}_{D'}^{-(t+1)}\otimes \mathcal{T}^{(t+1)m+j}\to
\mathcal{T}_{D'}^{-t}\otimes \mathcal{T}^{(t+1)m+j}\to
\mathcal{T}_{D'}^{-t}\otimes \mathcal{T}^{(t+1)m+j}\big|_{E'}\to 0$$
we obtain an exact sequence
$$\HH^1(\mathcal{T}_{D'}^{-t}\otimes \mathcal{T}^{(t+1)m+j})
\to \HH^1(\mathcal{T}_{D'}^{-t}\otimes \mathcal{T}^{(t+1)m+j}\big|_{E'}) \to 
\HH^2(f\mathcal{T}_{D'}^{-(t+1)}\otimes \mathcal{T}^{(t+1)m+j})
.$$
We have $\HH^1(\mathcal{T}_{D'}^{-t}\otimes \mathcal{T}^{(t+1)m+j})\cong
\HH^1(\mathcal{T}^{m+j})=0$.  Furthermore, $\HH^2(f\mathcal{T}_{D'}^{-(t+1)}\otimes \mathcal{T}^{(t+1)m+j})
\cong \HH^2(\mathcal{T}^j)$, and this is zero when $j = 0$ since $X$ is acyclic,
and zero for $j > 0$ by Hypothesis \ref{hypo:Lamp}.   Hence the right column of (\ref{eq:bigcohomdiagram}) is 
exact, which completes the proof of  Lemma \ref{lem:cohomdiagram}.
\end{proof}

\begin{subcor}
\label{cor:surfaceind}
Let $j$ be an integer with $0\leq j<m$. Then for any
positive integer $t$, there is an exact sequence of $kG$-modules
\begin{equation}
\label{eq:sessur}
0\to K \to
\HH^0(\mathcal{T}^{m+j}) \oplus
f\,\HH^0(\mathcal{T}^{tm+j}) \to \HH^0(\mathcal{T}^{(t+1)m+j}) \to
Q_{t,j}\to 0
\end{equation}
where $Q_{t,j}=\HH^0(\mathcal{T}^{(t+1)m+j}\big|_{tD'\cap E'})$ is a free $kG$-module,
and $K\cong f\,\HH^0(\mathcal{T}^j)$.
\end{subcor}

\begin{proof}
Since $E'\cap D'\cap B=\emptyset$, it follows that $tD'\cap E'$ is a $0$-dimensional subvariety of $X$ with  \'{e}tale Galois $G$-action. Thus $Q_{t,j}$ is a free $kG$-module
by Lemma \ref{lem:project}. 
 By Lemma 
\ref{lem:cohomdiagram}, there is an exact 
sequence of $kG$-modules
\begin{equation}
\label{eq:oyvey!}
\HH^0(\mathcal{T}_{D'}^{-t}\otimes \mathcal{T}^{(t+1)m+j}) \oplus
\HH^0(f\mathcal{T}_{D'}^{-1}\otimes\mathcal{T}^{(t+1)m+j}) \arrow{\varphi}{}
\HH^0(\mathcal{T}^{(t+1)m+j}) \to Q_{t,j}\to 0.
\end{equation}
The $G$-equivariant $\mathcal{O}_X$-module
isomorphism $\lambda:\mathcal{T}^m\to \mathcal{T}_{D'}$ induces isomorphisms
$$(\lambda^{-1})^{\otimes(t+1)}:  \mathcal{T}_{D'}^{t+1-s}\otimes \mathcal{T}^{j}\to 
\mathcal{T}_{D'}^{-s}\otimes\mathcal{T}^{(t+1)m+j}$$
for all $0\le j<m$ and for all $1\le s\le t$.
We obtain that
\begin{equation}
\label{eq:lost1}
\HH^0(\mathcal{T}_{D'}^{-t}\otimes \mathcal{T}^{(t+1)m+j}) = (\lambda^{-1})^{\otimes
(t+1)}\,\HH^0(\mathcal{T}_{D'}\otimes \mathcal{T}^{j}) = 
(\lambda^{-1})^{\otimes (t+1)} \,\lambda \,\HH^0(\mathcal{T}^{m+j})
\end{equation}
and 
\begin{equation}
\label{eq:lost2}
\HH^0(f\mathcal{T}_{D'}^{-1}\otimes\mathcal{T}^{(t+1)m+j})
= (\lambda^{-1})^{\otimes (t+1)}f\,\HH^0(\mathcal{T}_{D'}^{t}\otimes
\mathcal{T}^{j})\, =
(\lambda^{-1})^{\otimes (t+1)}\, \lambda^{\otimes t}\,
f\,\HH^0(\mathcal{T}^{tm+j})
\end{equation}
We now show that the intersection of  $\HH^0(\mathcal{T}_{D'}\otimes 
\mathcal{T}^{j})$ and $f\,\HH^0(\mathcal{T}_{D'}^{t}\otimes \mathcal{T}^{j})$ is
$f\,\HH^0(\mathcal{T}^j)$ when these cohomology groups are regarded as submodules of
$\HH^0(\mathcal{T}_{D'}^{t+1}\otimes\mathcal{T}^j)$.  By Hypothesis
\ref{hypo:Lamp}, $\mathcal{T}$ is isomorphic to 
$\mathcal{O}_X(T)$, where $T$ is effective.  We need to show that
the intersection of $\HH^0(\mathcal{T}_{D'}\otimes \mathcal{O}_X(jT))$
and $f\,\HH^0(\mathcal{T}_{D'}^{t}\otimes \mathcal{O}_X(jT))$ is
$f\,\HH^0(\mathcal{O}_X(jT))$, when these cohomology groups are regarded as submodules of
$\HH^0(\mathcal{T}_{D'}^{t+1}\otimes\mathcal{O}_X(jT))$. 
Let $h\in \HH^0(\mathcal{T}_{D'}^{t}\otimes 
\mathcal{O}_X(jT))$ with $fh\in \HH^0(\mathcal{T}_{D'}\otimes 
\mathcal{O}_X(jT))$.  We have $$0\leq 
\mathrm{div}_X(h)+tD'+jT\quad \mathrm{and}\quad 0\leq \mathrm{div}_X(fh)+D'+jT=\mathrm{div}_X(h)+E'
+jT.$$ Hence $\mathrm{div}_X(h)+jT\geq 0$, which means that $h\in\HH^0(\mathcal{O}_X(jT))$.
It follows that the kernel of $\varphi$ in (\ref{eq:oyvey!}) is $(\lambda^{-1})^{\otimes
(t+1)}f\,\HH^0(\mathcal{T}^j)$, so Corollary \ref{cor:surfaceind} now
follows from (\ref{eq:lost1}) and (\ref{eq:lost2}).
\end{proof}

\begin{subcor}
\label{cor:surface}
Let $j$ be an integer with $0\leq j<m=\#G$. Then for any
positive integer $t$, there is an isomorphism of $kG$-modules 
$$\HH^0(\mathcal{T}^{(t+1)m+j})\cong \frac{\HH^0(\mathcal{O}_X(D'+jT))^{t+1}}{K_{t,j}} \oplus 
P_{t,j}$$
where  $\HH^0(\mathcal{O}_X(D'+jT))^{t+1}$ is the direct sum of $t+1$ copies of
$\HH^0(\mathcal{O}_X(D'+jT))$, 
\begin{eqnarray*}
K_{t,j}&=&\HH^0(\mathcal{O}_X(jT))\cdot (f,-1,0,0,\ldots,0)\oplus 
\HH^0(\mathcal{O}_X(jT)) \cdot(0,f,-1,0\ldots,0) \\
&& \oplus \cdots \oplus \HH^0(\mathcal{O}_X(jT))
\cdot(0,\ldots,0,f,-1)
\end{eqnarray*}
and $P_{t,j}$ is a free $kG$-module.
\end{subcor}

\begin{proof}
Let $\lambda'={\lambda}^{-1}:\mathcal{O}_X(D')\to \mathcal{T}^m$ and let 
$\psi'={\psi}^{-1}:\mathcal{O}_X(T)\to \mathcal{T}$. Then for all $0\leq j<m$
and for all positive $t$, we
have a $G$-equivariant $\mathcal{O}_X$-module homomorphism
$${\lambda'}^{\otimes (t+1)}\otimes{\psi'}^{\otimes j}:\mathcal{O}_X((t+1)D'+jT)\to
\mathcal{T}^{(t+1)m+j}.$$ We prove by induction on $t$ that
\begin{equation}
\label{eq:whatsupdoc}
\HH^0(\mathcal{T}^{(t+1)m+j})={\lambda'}^{\otimes (t+1)}\otimes{\psi'}^{\otimes j}\left(
\frac{\HH^0(\mathcal{O}_X(D'+jT)) \oplus 
\cdots \oplus \HH^0(f^t\mathcal{O}_X(D'+jT))}{F_{t,j}} \oplus P_{t,j}\right)
\end{equation}
where 
\begin{eqnarray*}
F_{t,j}&=&\HH^0(\mathcal{O}_X(jT))\cdot (f,f\cdot(-1),0,0,\ldots,0)\oplus 
\HH^0(\mathcal{O}_X(jT)) \cdot(0,f\cdot f,f^2\cdot(-1),0\ldots,0) \\
&& \oplus \cdots \oplus \HH^0(\mathcal{O}_X(jT))
\cdot(0,\ldots,0,f^{t-1}\cdot f,f^t\cdot(-1)).
\end{eqnarray*}
Since $f\in k(Y)=k(X)^G$, there is a $kG$-module isomorphism
 $$\HH^0(f^l\mathcal{O}_X(D'+jT)) \cong  \HH^0(\mathcal{O}_X(D'+jT))$$ for all $l$.  Hence
Corollary \ref{cor:surface} will follow from (\ref{eq:whatsupdoc}).

When $t = 1$, (\ref{eq:whatsupdoc}) follows from the proof of Corollary \ref{cor:surfaceind}
with $P_{t,j} \cong Q_{1,j}$.  
For $t > 1$,
the proof of Corollary \ref{cor:surfaceind} shows that
\begin{eqnarray}
\label{eq:weahhh}
\HH^0({\mathcal{T}}^{(t+1)m+j})&=&{\lambda'}^{\otimes (t+1)}\otimes{\psi'}^{\otimes j}\left(\frac{
\HH^0(\mathcal{O}_X(D'+jT))\oplus 
f\,\HH^0(\mathcal{O}_X(tD'+jT))}{
\{\mu(f,f\cdot (-1))\,|\, \mu\in \HH^0(\mathcal{O}_X(jT))\}}\right)\oplus Q_{t,j}\nonumber\\
&=& {\lambda'}^{\otimes (t+1)}\otimes{\psi'}^{\otimes j}\left(\frac{
\HH^0(\mathcal{O}_X(D'+jT))\oplus 
(\psi^{\otimes j}\otimes\lambda^{\otimes t})\,f\,\HH^0(\mathcal{T}^{tm+j})}{
\{\mu(f,f\cdot (-1))\,|\, \mu\in \HH^0(\mathcal{O}_X(jT))\}}\right)\oplus Q_{t,j}.
\end{eqnarray}
By induction, 
\begin{equation}
\label{eq:inductstep}
(\psi^{\otimes j}\otimes\lambda^{\otimes t})\,f\,\HH^0(\mathcal{T}^{tm+j}) = 
f \left[\frac{\HH^0(\mathcal{O}_X(D'+jT)) \oplus \cdots 
\oplus \HH^0(f^{t-1}\mathcal{O}_X(D'+jT))}{F_{t-1,j}} \oplus 
P_{t-1,j}\right].
\end{equation}
Substituting (\ref{eq:inductstep}) into (\ref{eq:weahhh}) completes
the induction and the proof.
\end{proof}

The following Theorem follows from Corollary \ref{cor:surface}.

\begin{subthm}
\label{thm:surfacenec}
With the notation of  Hypothesis $\ref{hypo:Lamp}$,
$S(X,\mathcal{T})$ is an  indecomposably finite $kG$-module if and only if for each integer $0\leq j<m = \#G$ 
there are finitely many non-isomorphic finitely generated indecomposable $kG$-modules $U(j,1),\ldots, 
U(j,r_j)$ such that 
for all positive integers $t$ there exist positive integers $l_1(t,j),\ldots, l_{r_j}(t,j)$ with
$$\frac{\HH^0(\mathcal{O}_X(D'+jT))^{t+1}}{K_{t,j}}
\cong U(j,1)^{l_1(t,j)}\oplus\cdots\oplus U(j,r_j)^{l_{r_j}(t,j)}$$
as $kG$-modules, 
where $K_{t,j}$ is as in Corollary $\ref{cor:surface}$.
\end{subthm}

Our goal now is to deduce from this Theorem the following result,
which proves part (ii) of Theorem \ref{thm:main}.

\begin{subthm}
\label{thm:surface1}
With the notation of  Hypothesis $\ref{hypo:Lamp}$,
 $S(X,\mathcal{T}^{\#G})$ is an indecomposably finite $kG$-module
having a polynomial description of degree $d=2$.
\end{subthm}

\begin{subrem}
\label{rem:wedontknow}
We do not know whether or not $S(X,\mathcal{T})$ is indecomposably finite.
Theorem \ref{thm:surfacenec} shows that this only depends on the $kG$-module structure of $\HH^0(\mathcal{O}_X(
D'+jT))$  and on how the elements $f$ and $1$ lie in $\HH^0(\mathcal{O}_X(
D'+jT))$ for $0\leq j< m = \#G$.
\end{subrem}

To prove Theorem \ref{thm:surface1}, we need 
the following result.

\begin{sublem}
\label{lem:biteme}
Suppose  that $Z$ is a smooth projective geometrically integral acyclic variety over a field $L$ of
characteristic $p$ on which a finite $p$-group
$G_p$ acts over $L$.   Then there is a closed point $z_0 \in Z$ which is fixed by $G_p$
as a closed subscheme of $Z$.  If $L$ is finite, then $z_0$ can be chosen to have residue field $L$.
\end{sublem}

\begin{proof}  Suppose first that $L = \mathbb{F}_q$ for some power $q$ of $p$, and let $\overline{L}$
be an algebraic closure of $L$.  The zeta function of $Z$ is
$$\zeta(Z,t) = \prod_{x \in Z^0} (1 - t^{\mathrm{deg}(x)})^{-1}$$
where the product is over the set $Z^0$ of closed points of $Z$,
and $\mathrm{deg}(x)$ is the degree of $x$ over $L$.  Define $\overline{Z} = \overline{L} \otimes_L Z$ and let $\mathbb{F}_q$ be the constant sheaf defined by $L$ for the \'etale topology
of $Z$.  Define $F^*$ to be the geometric Frobenius automorphism of $\HH^n_{et}(\overline{Z},\mathbb{F}_q)$
for all $n$.  By a result of Katz (see \cite[XXII, 3.1]{SGA7} and \cite[Thm 2.2]{SGA4.5}),
there is a congruence
\begin{equation}
\label{eq:zappa}
\zeta(Z,t) \equiv \prod_{n \ge 0} \mathrm{det}_{\mathbb{F}_q}(1 - F^* t | \HH^n_{et}(\overline{Z},\mathbb{F}_q))^{(-1)^{n+1}}\mod p
\end{equation}
where on the right, one takes determinants of endomorphisms of $\mathbb{F}_q$-vector spaces.
As in \cite[2.18(c)]{Milne}, we have an Artin Schreier sequence of sheaves on $Z_{et}$
\begin{equation}
\label{eq:ashr}
0 \arrow{}{} \mathbb{F}_q \arrow{}{} \mathbb{G}_a \arrow{f - 1}{} \mathbb{G}_a \arrow{}{} 0
\end{equation}
where $f$ is the map $a \to a^q$.  By \cite[Prop. 3.7, Them. 3.9]{Milne}, $\HH^n_{et}(\overline{Z},\mathbb{G}_a) = \HH^n(\overline{Z},\mathcal{O}_{\overline{Z}})$ for all $n$.  Since $Z$ is acyclic and geometrically
integral over the field $L = \mathbb{F}_q$, this implies
$\HH^n_{et}(\overline{Z},\mathbb{G}_a)$ is $\overline{\mathbb{F}}_q$ (resp. $0$) if $n = 0$ (resp. 
$n > 0$).
Taking cohomology of \ref{eq:ashr} over $\overline{Z}$ now shows that $\HH^n_{et}(\overline{Z},\mathbb{F}_q)$ 
is $\mathbb{F}_q$ (resp. $0$) if $n = 0 $ (resp. if $n > 0$).  Thus (\ref{eq:zappa}) gives 
the congruence
$$\zeta(Z,t) = \prod_{x \in Z^0} (1 - t^{\mathrm{deg}(x)})^{-1} \equiv (1 - t)^{-1} = \sum_{m \ge 0} t^m 
\mod p.$$
This implies in particular that $\# Z(L) = \#\{x \in Z^0: \mathrm{deg}(x) = 1\} \equiv 1 $ mod $p$.  Since $G_p$ is a $p$-group
and permutes $Z(L)$, it has to fix a point of $Z(L)$.  

If $L$ is of finite transcendence degree over $\mathbb{F}_p$,
one chooses a quasi-projective variety $\mathcal{V}$ over $\mathbb{F}_p$ with
function field $L$, and an integral model $\mathcal{Z}$ of $Z$
over $\mathcal{V}$ which is flat over $\mathcal{V}$.  We can find
a dense open subset $\mathcal{V}_0$ of closed points $v \in \mathcal{V}$ such that
the fiber $\mathcal{Z}_v = \mathrm{Spec}(k(v)) \times \mathcal{Z}$ is
a smooth projective geometrically integral acyclic variety over the finite residue field $k(v)$ of $v$.   Since
$G_p$ acts on $\mathcal{Z}_v$, the above argument shows that $\mathcal{Z}_v^{G_p}$
is non-empty.  Hence the fixed point variety $\mathcal{Z}^{G_p}$
projects to a subset of $\mathcal{V}$ which contains a dense
set of closed points.  However, from \cite[p. 94]{hart}, we know that
$\mathcal{Z}^{G_p}$ is closed
in $\mathcal{Z}$ and has constructible image in $\mathcal{V}$, and that
a constructible subset of $\mathcal{V}$ is dense if and only
if it contains the generic point.  We conclude that $\mathcal{Z}^{G_p}$
has image containing the generic point, so $Z^{G_p}$ is non-empty.

Finally, for arbitrary $L$ of characteristic $p$, we use the
fact that $G_p$ is finite to deduce that $Z$ and the $G_p$-action
on $Z$ are defined over a field of finite transcendence
degree over $\mathbb{F}_p$.  Hence we can reduce to the previous case.
\end{proof}

\medskip

\noindent
\textit{Proof of Theorem $\ref{thm:surface1}$.}
Setting $j = 0 $ in Corollary \ref{cor:surface} shows that for any positive integer $t$
\begin{equation}
\label{eq:olala}
\HH^0(\mathcal{T}^{(t+1)m}) \cong \frac{\HH^0(\mathcal{T}_{D'})^{t+1}}
{K_{t}} \oplus P_{t}
\end{equation}
where 
\begin{eqnarray*}
K_{t}&=&k\cdot (f,-1,0,0,\ldots,0)\oplus k\cdot
(0,f,-1,0\ldots,0)  \oplus \cdots \oplus k\cdot (0,\ldots,0, f,-1)
\end{eqnarray*}
and $P_{t}$ is a free $kG$-module.

By Lemma \ref{lem:biteme}, there is a closed point $x_0 \in X$ fixed by a Sylow $p$-subgroup $G_p$ of $G$.  In Hypothesis \ref{hypo:Lamp} we let $A$ be the $G$-orbit of $x_0$.  
Thus $A \cap D' = A \cap E' = \emptyset$ by Hypothesis \ref{hypo:Lamp}, so $f \in \mathcal{O}_{X,a}^*$ for $a \in A$ because   $\mathrm{div}_X(f)=E'-D'$.  Let $s(a) = 
(\mathcal{T}_{D'})_a/ \left(\mathfrak{m}_{X,a}\cdot 
(\mathcal{T}_{D'})_a\right) \cong k$ for $a \in A$.  Let 
$H$ be the inertia group of $x_0$, so that $G_p \subset H$.  Consider
the natural $kG$-module homomorphism
\begin{equation}
\label{eq:into}
kf \to \bigoplus_{a \in A} s(a) = kG \otimes_{kH} s(x_0)
\end{equation}
induced by taking the image of $f$ in the stalks of $\mathcal{T}_{D'}$.  Since
the image of $f$ in $s(x_0) \cong k$ is non-trivial
and since $kf$ is $(G,H)$-injective by \cite[Thm. 19.2 and Prop. 19.5]{CR},
we find 
that (\ref{eq:into}) is split by a $kG$-module homomorphism
$\tau:  \bigoplus_{a \in A} s(a) \to kf$.  Let $\nu:\HH^0(\mathcal{T}_{D'}) \to kf$
be the composition of the 
natural $kG$-module homomorphism
$\HH^0(\mathcal{T}_{D'})\to \bigoplus_{a \in A} s(a)$ with $\tau$. Define $U$ to be the kernel of $\nu$. Then $k\cdot f$ is a $kG$-complement 
of $U$ in $\HH^0(\mathcal{T}_{D'})$, and $\HH^0(\mathcal{T}_{D'}) = U\oplus k f$ as 
$kG$-modules.

Let $U_t=U^{t} \oplus \HH^0(\mathcal{T}_{D'})$, considered as a $kG$-submodule
of $\HH^0(\mathcal{T}_{D'})^{t+1}$, and suppose that 
$$(f_1,\ldots,f_t,f_{t+1})\in K_t\cap U_t.$$ 
Then there exist elements $a_1,\cdots,a_t\in K$ such that
$(f_1,\cdots,f_t,f_{t+1})=(a_1 f, a_2 f - a_1,\ldots, a_t f - a_{t-1}, - a_t) \in U_t$. Therefore 
$a_1 f, a_2 f - a_1,\ldots, a_t f - a_{t-1}$ are elements of $U$. Since $U \cap k f =0$,
it follows that $a_1=0=a_2=\cdots = a_{t-1}=a_t$. Hence $K_t\cap U_t=0$. On the other hand,
since the $k$-dimension of $\HH^0(\mathcal{T}_{D'})$ is $\mathrm{dim}_k (U) +1$,
it follows that $\mathrm{dim}_k(K_t) + \mathrm{dim}_k(U_t) = (t+1)\,\mathrm{dim}_k
(\HH^0(\mathcal{T}_{D'}))$. Hence  
$$\frac{\HH^0(\mathcal{T}_{D'})^{t+1}}{K_t}\cong U_t = U^{t} \oplus 
\HH^0(\mathcal{T}_{D'}) $$
as $kG$-modules.  This and (\ref{eq:olala}) give a $kG$-module isomorphism
\begin{equation}
\label{eq:olala2}
\HH^0(\mathcal{T}^{(t+1)m}) \cong U^{t} \oplus 
\HH^0(\mathcal{T}_{D'}) \oplus P_{t}
\end{equation}
where $P_t$ is a free $kG$-module. Since $\mathrm{dim}_k( \HH^0(\mathcal{T}^{(t+1)m}))$ grows as a quadratic 
polynomial in $t$ for large $t$, this implies
$S(X,\mathcal{T}^m)$
is indecomposably finite and has a polynomial description of degree $d=2$.
 \hfill \BBox

\subsection{Projective space}
\label{ss:proj}

In this subsection we prove part (iii) of Theorem \ref{thm:main}. 
Let $X=\PP_k^d=\mathrm{Proj}\;k[z_0,z_1,\cdots,z_d]$.    
We assume that $G$ acts on $\PP_k^d$ over $k$ and that this action is 
faithful and generically free such that
$\mathcal{L}=\mathcal{O}(1)$ is a (very ample) $G$-equivariant line bundle on $\PP_k^d$.
By Lemma \ref{lem:biggerfield}, we
can assume $k$ is an algebraically closed field of characteristic $p > 0$. 
After a linear change of the homogeneous coordinates $z_0,\ldots,z_d$, we may
assume that relative to these coordinates the action of elements in a Sylow $p$-subgroup
$G_p$ of $G$ on 
$\HH^0(\mathcal{L})=k z_0+\cdots+ k z_d$ is given by upper triangular unipotent  matrices. 
Define $T$ to be the hyperplane $z_0=0$, and $Q$ to be the point $(0:\cdots:0:1)$ 
in $\PP_k^d$ over $k$. Then both $T$ and $Q$ are $G_p$-stable.
Let $Y=\PP_k^d/G$, define $\pi:\PP_k^d\to Y$ to be the quotient morphism, and let $b$ (resp. $B$)
be the branch locus in $Y$ (resp. the ramification locus in $\PP_k^d$) of $\pi$.

Our goal is to construct a complex with $d+2$ terms which is an exact Koszul resolution,
and then to study its splitting behavior.

\begin{subdfn}
\label{def:koszul}
By Lemma \ref{lem:linebundle}, there exist effective Weil divisors 
$D_0,\ldots,D_{d-1}$ 
on $Y$ such that when $D_j'=\pi^*(D_j)$ the following is true:
\begin{enumerate}
\item[a.] for all $j$, $D_j'$ does not contain any point in the $G$-orbit of $Q=(0:\cdots:0:1)$;
\item[b.] there is a $G$-equivariant $\mathcal{O}_{\PP_k^d}$-module isomorphism 
$$\lambda:\quad  \mathcal{L}^{m}\to \mathcal{O}_{\PP_k^d}(D_0');$$
where $m= m_2^2 (\#G)$ and $m_2$ is the maximal prime to $p$ divisor of $\# G$.
\item[c.] the divisors $D_0,\ldots,D_{d-1}$ are all linearly equivalent on $Y$;
\item[d.] for each subset $J\subset \{0,\ldots,d-1\}$ and each point 
$x \in \bigcap_{j \in J} D'_j$, the local equations of the $D'_j$ in 
$\mathcal{O}_{{\PP_k^d},x}$ for $j \in J$ form
a regular sequence; 
\item[e.] the intersection $D_0'\cap \cdots \cap D_{d-1}'\cap B=\emptyset$.
\end{enumerate}
Let $D=D_0$, $D'=D_0'$ and let $\mathcal{L}_{D'}=\mathcal{O}_{\PP_k^d}(D')$.
Define $f_0=1$, and let  $f_1,f_2,\ldots, f_{d-1}$ be in $k(Y)$ with $\mathrm{div}_Y(f_j)=D_j-D_0$. 
The Koszul resolution for $\mathcal{O}_{D_0'
\cap D_1'\cap\cdots \cap D_{d-1}'}$ has
the following form
\begin{equation}
\label{eq:koszul1}
0\to K_d \xrightarrow{\rho_d} K_{d-1} \xrightarrow{\rho_{d-1}} \cdots 
\xrightarrow{\rho_2} K_1 \xrightarrow{\rho_1} K_0
\xrightarrow{\rho_0}  \mathcal{O}_{D_0'\cap D_1'\cap\cdots \cap D_{d-1}'}\to 0
\end{equation}
where $K_0=\mathcal{O}_{\PP_k^d}$, 
$$K_1=\mathcal{O}_{\PP_k^d}(-D_0')\oplus \mathcal{O}_{\PP_k^d}(-D_1')\oplus \cdots \oplus
\mathcal{O}_{\PP_k^d}(-D_{d-1}') = \bigoplus_{0\leq i\leq d-1} f_i\mathcal{L}_{D'}^{-1},$$
and, for $2\leq r\leq d$,
$$K_r=\bigwedge\!{}^r K_1 = \bigoplus_{0\leq i_1<i_2<\cdots <i_r\leq d-1} f_{i_1}f_{i_2}\cdots  f_{i_r}
\mathcal{L}_{D'}^{-r}.$$
The morphisms $\rho_1$ and $\rho_0$ are the obvious ones, and for $2\leq r\leq d$,
 $\rho_r$ is defined by ``contraction,'' namely
 $$\rho_r(t_1\wedge \cdots \wedge t_r) = \sum_{j=1}^r (-1)^{j-1} \rho_1(t_j) \;(t_1\wedge
 \cdots \wedge \widehat{t_j} \wedge\cdots \wedge t_r).$$
Because of property (d) above the Koszul resolution is exact (see \cite[Thm. 43, p. 135]{Matsumura}).
Tensoring the Koszul resolution with any positive power of 
$\mathcal{L}$ results in an exact sequence. 
Since $\mathcal{L}$ is very ample, there exists an integer $\mu_0\geq d$ such that 
for all $0\leq r\leq d$, $\HH^1(
\mathrm{ker}(\rho_r)\otimes\mathcal{L}^n)=0$ for all $n\geq m\mu_0$. 
Let now $j$ be an integer with $0\leq j < m$, and let $t\geq \mu_0$ be 
an integer.   
The isomorphism $\lambda: \mathcal{L}^{m}\to \mathcal{O}_{\PP_k^d}(D_0')=
\mathcal{O}_{\PP_k^d}(D')$
gives an isomorphism $\lambda^{\otimes t}:\mathcal{L}^{tm+j} \to \mathcal{L}_{D'}^t\otimes\mathcal{L}^j$.
By tensoring the Koszul resolution (\ref{eq:koszul1}) with 
$\mathcal{L}_{D'}^t\otimes\mathcal{L}^j=\lambda^{\otimes t}\,\mathcal{L}^{tm+j}$ and 
taking cohomology on $\PP_k^d$ of the resulting complex, we obtain
the following exact Koszul complex of $kG$-modules, denoted by
$C_\bullet=C_\bullet(t,j)$,
\begin{equation}
\label{eq:koszul}
0\to C_d \xrightarrow{\tau_d} C_{d-1} \xrightarrow{\tau_{d-1}} \cdots \xrightarrow{\tau_2} C_1\xrightarrow{\tau_1} C_0\xrightarrow{\tau_0}  Q_{t,j}\to 0
\end{equation}
where $Q_{t,j}=\HH^0(\mathcal{L}_{D'}^{t}\otimes \mathcal{L}^j\big|_{D_0'
\cap D_1'\cap\cdots \cap D_{d-1}'})$, 
$C_0=\HH^0(\mathcal{L}_{D'}^{t}\otimes\mathcal{L}^j)$, and, for $1\leq r\leq d$,
\begin{equation}
\label{eq:crdef}
C_r= \bigoplus_{0\leq i_1<i_2<\cdots <i_r\leq d-1} f_{i_1}f_{i_2}\cdots  f_{i_r}\;
\HH^0( \mathcal{L}_{D'}^{t-r}\otimes\mathcal{L}^j).
\end{equation}
The morphism $\tau_0$ is the canonical projection, and, for
$1\leq r\leq d$, $\tau_r$ is defined by
$$\tau_r(f_{i_1}f_{i_2}\cdots  f_{i_r}\cdot h) = \sum_{j=1}^r (-1)^{j-1} f_{i_j}\;(f_{i_1}\cdots 
\widehat{f_{i_j}}\cdots f_{i_r} \cdot h)$$
for $h\in \HH^0( \mathcal{L}_{D'}^{t-r}\otimes\mathcal{L}^j)$.

Note that since  $D_0' \cap D_1'\cap\cdots \cap D_{d-1}'$ is a $0$-dimensional subvariety 
of $\PP_k^d$ with  \'{e}tale Galois $G$-action, $Q_{t,j}$ is a free $kG$-module by Lemma \ref{lem:project}. 
\end{subdfn}

\begin{sublem}
\label{lem:verynice}
Let $\mathcal{L}_{D'}=\mathcal{O}_{\PP_k^d}(D')$.  For $x$ a closed point of $X$, 
$t\geq \mu_0$, $0\leq j<m$ and $1\leq r\leq d$ let
$\mathfrak{m}_{\PP_k^d,x}$ be the maximal ideal of $\mathcal{O}_{\PP_k^d,x}$, and define
$$s(t,j,r,x) = \frac{(\mathcal{L}_{D'}^{t-r+1}\otimes \mathcal{L}^{j})_{x}}
{\mathfrak{m}_{\PP_k^d,x}^{(t-r)m+j+1}\cdot(\mathcal{L}_{D'}^{t-r+1}\otimes
\mathcal{L}^{j})_{x}}.$$
Let $H$ be the inertia group of the point
$Q= (0:\cdots:0:1)$, so that $G_p \subset H$.   Then for all $0\leq i \leq d-1$,
the natural $kH$-module homomorphism
$$ \sigma_{t-r,j,i,Q}: \quad f_i\;\HH^0(\mathcal{L}_{D'}^{t-r}\otimes \mathcal{L}^{j})
\quad\to\quad  s(t,j,r,Q)$$
is an isomorphism.
\end{sublem}

\begin{proof}
Since $\sigma_{t-r,j,i,Q}$ is a $kH$-module homomorphism, it is enough to show it is 
$k$-vector space isomorphism.   Because $D'$ is linearly equivalent to
the divisor of $z_d^m$, we can reduce to showing $\sigma_{t-r,j,i,Q}$ is a $k$-vector
space isomorphism when $D'$ is the hypersurface $z_d^m=0$. Since $Q$ does not
lie on $z_d^m = 0$,  the constant $1$ is a local equation
for $D'$ at $Q$, and $(z_0/z_d)$ is a local equation for $T$ at $Q$. We obtain
\begin{eqnarray*}
\frac{(\mathcal{L}_{D'}^{t-r+1}\otimes \mathcal{L}^{j})_{Q}} {\mathfrak{m}_{\PP_k^d,Q}^{
(t-r)m+j+1}\cdot(\mathcal{L}_{D'}^{t-r+1}\otimes\mathcal{L}^{j})_{Q}}
&=&\frac{ \left(\frac{z_d}{z_0}\right)^{j}\cdot \mathcal{O}_{\PP_k^d,Q}}
{  \left(\frac{z_d}{z_0}\right)^{j}\cdot \mathfrak{m}_{\PP_k^d,Q}^{(t-r)m+j+1}}
\;=\;\frac{\left(\frac{z_d}{z_0}\right)^{j}f_i\cdot \mathcal{O}_{\PP_k^d,Q}}{\left(\frac{z_d}{z_0}\right)^{j}f_i
\cdot\mathfrak{m}_{\PP_k^d,Q}^{(t-r)m+j+1}}\\
&=& \left(\frac{z_d}{z_0}\right)^{j}f_i\cdot \frac{\mathcal{O}_{\PP_k^d,Q}}{\mathfrak{m}_{\PP_k^d,
Q}^{(t-r)m+j+1}} 
\end{eqnarray*}
where 
the second equation follows since  $f_i$ is a unit in $\mathcal{O}_{\PP_k^d,Q}$.
Then 
$$\frac{\mathcal{O}_{\PP_k^d,Q}}{\mathfrak{m}_{\PP_k^d,Q}^{(t-r)m+j+1}}=
\frac{k[\frac{z_0}{z_d},\ldots,\frac{z_{d-1}}{z_d}]}{(\frac{z_0}{z_d},\ldots,\frac{z_{d-1}}{z_d})^{(t-r)
m+j+1}}$$
and the latter module can be identified with the space of all polynomials in $\displaystyle
\frac{z_0}{z_d},\ldots, \frac{z_{d-1}}{z_d}$ of degree $\leq (t-r)m+j$.
Now let $h\in\HH^0(\mathcal{L}_{D'}^{t-r}\otimes\mathcal{L}^j)$. This means that
$\mathrm{div}_X(h)+(t-r)D'+jT\geq 0$,
i.e. $h=g/(z_d^{(t-r)m}\,z_0^j)$ where $g$ is a homogeneous polynomial
in $k[z_0,z_1,\ldots,z_d]$ of degree $(t-r)m+j$. Then $f_i\,h=\displaystyle 
\left(\frac{z_d}{z_0}\right)^jf_i \cdot \frac{g}{z_d^{(t-r)m+j}}$, and $\displaystyle 
\frac{g}{z_d^{(t-r)m+j}}$ can be identified
with a polynomial in $\displaystyle \frac{z_0}{z_d},\ldots, \frac{z_{d-1}}{z_d}$ of degree $\leq
(t-r)m+j$. Also, as $h$ runs over all the elements in $\HH^0(\mathcal{L}_{D'}^{t-
r}\otimes\mathcal{L}^j)$, all polynomials in $\displaystyle \frac{z_0}{z_d},\ldots,
 \frac{z_{d-1}}{z_d}$ of degree $\leq (t-r)m+j$ occur.
Hence the morphism $\sigma_{t-r,j,i,Q}$ from
$f_i\;\HH^0(\mathcal{L}_{D'}^{t-r}\otimes \mathcal{L}^{j})$ to $\displaystyle
\frac{(\mathcal{L}_{D'}^{t-r+1}\otimes \mathcal{L}^{j})_{Q}} {\mathfrak{m}_{\PP_k^d,Q}^{(t-r)m+j+1}\cdot(\mathcal{L}_{D'}^{t-r+1}\otimes\mathcal{L}^{j})_{Q}}$ is an isomorphism of 
$k$-vector spaces.
\end{proof}

We have the following Corollary from Lemma \ref{lem:verynice}
and the fact that every finitely generated $kG$-module is $(G,H)$-injective
when $G_p\subset H$ \cite[Thm. 19.2 and Prop. 19.5]{CR}.

\begin{subcor}
\label{cor:verynice}
Suppose $\mathcal{L}_{D'}$,  $t$, $j$, $r$, $Q$, $H$ and $i$ are as 
in Lemma $\ref{lem:verynice}$.  Let $\omega(Q)$ be the $G$-orbit of $Q$.  The
$\sigma_{t-r,j,i,Q'}$ over $Q' \in \omega(Q)$ give a split
injection 
$$ \sigma_{t-r,j,i}: \quad f_i\;\HH^0(\mathcal{L}_{D'}^{t-r}\otimes \mathcal{L}^{j})
\quad\to\quad  \bigoplus_{Q' \in \omega(Q)} s(t,j,r,Q') = s(t,j,r)$$
of $kG$-modules.  Let $\tau: s(t,j,r) \to f_i\;\HH^0(\mathcal{L}_{D'}^{t-r}\otimes \mathcal{L}^{j})$
be a splitting homomorphism.  
The natural restriction homomorphisms
$\HH^0(\mathcal{L}_{D'}^{t-r+1}\otimes \mathcal{L}^{j}) \to
(\mathcal{L}_{D'}^{t-r+1}\otimes \mathcal{L}^{j})_{Q'}$
associated to $Q'\in\omega(Q)$ give a $kG$-module homomorphism
\begin{equation}
\label{eq:oyoyoy}
\HH^0(\mathcal{L}_{D'}^{t-r+1}\otimes \mathcal{L}^{j}) \quad\to\quad
s(t,j,r).
\end{equation}
Let
$$\widehat{\sigma}_{t-r,j}:\quad \HH^0(\mathcal{L}_{D'}^{t-r+1}\otimes \mathcal{L}^{j})
\quad \to\quad f_i\;\HH^0(\mathcal{L}_{D'}^{t-r}\otimes \mathcal{L}^{j})$$
be the composition of $(\ref{eq:oyoyoy})$ with $\tau$, and let $U_{t-r+1,j}=
\mathrm{ker}(\widehat{\sigma}_{t-r,j})$. 
Then we have an internal direct sum as $kG$-modules
$$\HH^0(\mathcal{L}_{D'}^{t-r+1}\otimes \mathcal{L}^{j}) =
U_{t-r+1,j} \oplus f_i\; \HH^0(\mathcal{L}_{D'}^{t-r}\otimes 
\mathcal{L}^{j}).$$
\end{subcor}

We now prove the following result for $\PP_k^2$ and $\PP_k^3$, which  also proves part (iii) of Theorem \ref{thm:main}.

\begin{subthm}
\label{thm:p23}
Let $d\in\{2,3\}$. Then the Koszul complex $C_{\bullet}(t,j)$ in $(\ref{eq:koszul})$ splits  
completely as a complex of $kG$-modules for all $0\leq j<m$ and for all $t\geq \mu_0$.
Moreover, the  polynomial ring $k[z_0,z_1,\ldots,z_d]$
is an indecomposably finite $kG$-module with a polynomial description of degree $d$.
\end{subthm}

\begin{subrem}
\label{rem:better}
Suppose $k$ is transcendental over $\mathbb{F}_p$. Then already in case of the projective
line $\PP_k^1$ over $k$, there may exist finite group actions which cannot be realized over any
algebraic extension of $\mathbb{F}_p$. For example,  suppose the characteristic of $k$ is $2$ and
$G=\mathbb{Z}/2\times\mathbb{Z}/2=\langle x,y\rangle$.  For each $c \in k$ one has an indecomposable representation $\rho_c$ defined by
$$x\mapsto \left[\begin{array}{cc}1&1\\0&1\end{array}\right],\qquad
y\mapsto \left[\begin{array}{cc}1&c\\0&1\end{array}\right]$$
and $\rho_c$ and $\rho_{c'}$ are isomorphic if and only if $c = c'$.  For this reason, Theorem \ref{thm:p23} does not follow from  Karagueuzian's and Symonds'
results in \cite{symkar2,sym1} in case $d=3$. However, the case $d=2$ was proved 
for an arbitrary field $k$ of characteristic $p$
in \cite{symkar2.5}.
\end{subrem}

To prove Theorem \ref{thm:p23} we will use the following Lemma.

\begin{sublem}
\label{lem:splitit} Suppose 
\begin{equation}
\label{eq:Mseq}
0 \to M_n \xrightarrow{\delta_n} M_{n-1} \xrightarrow{\delta_{n-1}} \cdots 
\xrightarrow{\delta_2} M_1 
\xrightarrow{\delta_{1}} M_0\to 0
\end{equation}
is an exact sequence of $kG$-modules for which $M_0$ is projective 
and for which there is an internal direct sum decomposition $M_i = 
\bigoplus_{l = 1}^{s(i)} M_{i,l}$ for each $i > 1$ having the following property.  
Let $\pi_{i,l}:M_i \to M_{i,l}$ be the $l^{th}$
projection and let $\iota_{i,l}:M_{i,l} \to M_i$ be the $l^{th}$ inclusion.  Then
if $i > 2$, there is an injective map $y_i:\{1,\cdots,s(i)\} \to \{1,\cdots,s(i-1)\}$ such that for all $1\le l\le s(i)$
\begin{equation}
\label{eq:inj}
\lambda_{i,l} = \pi_{i-1,y_i(l)} \circ \delta_i \circ \iota_{i,l} : \;
M_{i,l}\; \to \;
M_{i-1,y_i(l)}
\end{equation}
is a split injection.
Under these conditions, $(\ref{eq:Mseq})$ splits completely as a 
sequence of $kG$-modules.
\end{sublem}

\begin{proof}
The result is clear for $n\le 2$ and if the pair $(n,s(n))=(3,0)$.
Hence we use induction on $(n,s(n))$. 
Let $W_{n-1,y_n(s(n))}$ be a 
complement in $M_{n-1,y_n(s_n))}$ for the image of $\lambda_{n,s(n)}$ in 
(\ref{eq:inj}). Define $M'_{n-1}$ to be the inner direct sum
$$M'_{n-1} = \bigoplus_{\genfrac{}{}{0pt}{}{1\le i\le s(n-1)}{i\neq y_n(s(n))}}M_{n-1,i}\;\oplus\; W_{n-1,y_n(s(n))},$$
and let  $\iota'_{n-1}: M'_{n-1}\to M_{n-1}$ be the natural inclusion.
Using that $\lambda_{n,s(n)}$ from (\ref{eq:inj}) is a split injection, it follows 
that there is a short exact sequence of $kG$-modules
$$0 \to \delta_n\, M_{n,s(n)} \hookrightarrow M_{n-1}
\xrightarrow{\zeta'_{n-1}} M'_{n-1}\to 0$$
such that $\zeta'_{n-1}\circ \iota'_{n-1}$ is the identity on $M'_{n-1}$.
Since $M_{n-1}$ is the internal direct sum
$M_{n-1}=\delta_n\, M_{n,s(n)} \, \oplus \, \iota'_{n-1}\, M'_{n-1}$,
it follows that $\delta_{n-1}$ and $\delta'_{n-1}=\delta_{n-1}\circ
\iota'_{n-1}$ have the same image, $N_{n-2}$, in $M_{n-2}$. 
It also follows that $\delta_{n-1} = \delta'_{n-1}\circ\zeta'_{n-1}$.
We obtain a commutative diagram 
\begin{equation}
\label{eq:crazy}
\xymatrix @-.8pc {
0 \ar[r] & {\displaystyle\bigoplus_{i=1}^{s(n)-1} M_{n,i}} 
 \oplus M_{n,s(n)} \ar[r]^(.7){\delta_{n}} \ar@<2ex>[d]_{1\;\oplus\; 0} &
M_{n-1} \ar[d]^{\zeta'_{n-1}} \ar[r]^{\delta_{n-1}} & N_{n-2} \ar[r] 
\ar@{=}[d]& 0\\
0 \ar[r] & {\displaystyle\bigoplus_{i=1}^{s(n)-1} M_{n,i}}  \oplus {\;\;\;\; 0 \;\;\;\;}
\ar[r]^(.7){\delta_n} & M'_{n-1} \ar[r]^{\delta'_{n-1}} & N_{n-2} \ar[r] & 0
}
\end{equation}
By induction, the bottom row splits, i.e. there is a $kG$-module
homomorphism $\alpha'_{n-1}:N_{n-2}\to M'_{n-1}$ such that
$\delta'_{n-1}\circ\alpha'_{n-1}$ is the identity on $N_{n-2}$.
Then $\alpha_{n-1}=\iota'_{n-1}\circ\alpha'_{n-1}:N_{n-2}
\to M_{n-1}$ provides a $kG$-module splitting of the top row.
By induction, we also get that the exact sequence
$$0 \to {\displaystyle\bigoplus_{i=1}^{s(n)-1} M_{n,i}}  \xrightarrow{\delta'_n} M'_{n-1} \xrightarrow{\delta'_{n-1}} M_{n-2} \xrightarrow{\delta_{n-2}}
\cdots \xrightarrow{\delta_2} M_1 \xrightarrow{\delta_{1}} M_0\to 0$$
splits completely, where $\delta'_n$ is the restriction of $\delta_n$ to the 
submodule $\bigoplus_{i=1}^{s(n)-1} M_{n,i}$ of $M_n$. 
Hence the splitting of the top row in (\ref{eq:crazy})
implies that the sequence (\ref{eq:Mseq}) also splits completely.
\end{proof}

\medskip

\noindent
\textit{Proof of Theorem $\ref{thm:p23}$.}
Let $d=2$ or $3$.
We first want to use Lemma \ref{lem:splitit} to show 
that the Koszul complex $C_{\bullet}(t,j)$ in $(\ref{eq:koszul})$ splits  
completely as a complex of $kG$-modules for all $0\leq j<m$ and 
for all $t\geq \mu_0$. 
When $d=2$, the Koszul complex $C_{\bullet}(t,j)$ is a sequence
of the form (\ref{eq:Mseq}) with $n=3$, $s(3)=1$ and $s(2)=2$.
When $d=3$, $C_{\bullet}(t,j)$ is a sequence
of the form (\ref{eq:Mseq}) with $n=4$, $s(4)=1$ and $s(3)=3=s(2)$.
It follows from Corollary \ref{cor:verynice} that for $d=2,3$ the Koszul
complex $C_{\bullet}(t,j)$ satisfies the conditions of Lemma \ref{lem:splitit}.
This implies that $C_{\bullet}(t,j)$ splits  
completely as a complex of $kG$-modules for all $0\leq j<m$ and 
for all $t\geq \mu_0$. 

Hence for all $t\geq \mu_0$ and for all $0\leq j<m$, 
we have from (\ref{eq:koszul}) and (\ref{eq:crdef}) the following equality in $G_0^{\oplus}(kG)$:
$$[Q_{t,j}] 
= \sum_{r=0}^d (-1)^{r} [C_r]\nonumber\\
=  \sum_{r=0}^d (-1)^{r} {d \choose r}\cdot [\HH^0( \mathcal{L}^{m(t-r)+j})]$$
where $Q_{t,j}$ is a free $kG$-module.
Since $\mathrm{dim}_k( \HH^0(\mathcal{L}^{mt+j}))$ grows as a 
polynomial of degree $d$ in $t$ for large $t$, this implies
$k[z_0,\ldots,z_d]=\bigoplus_{n\geq 0}\HH^0(\mathcal{L}^n)$
is indecomposably finite and has a polynomial description of degree $d$.
\hfill \BBox


\end{document}